\documentclass[12pt]{amsart}
\usepackage[top=72pt,bottom=96pt,left=72pt,right=72pt]{geometry}
\usepackage{amsmath}
\usepackage{mathtools}
\usepackage{amssymb}
\usepackage{bbm}
\usepackage[all]{xy}
\usepackage{extarrows}

\def\G{\mathcal{G}}
\def\C{\mathbb{C}}

\def\N{\mathbb{N}}
\def\R{\mathbb{R}}
\def\F{\mathfrak{F}}
\def\f{\mathfrak{f}}

\def\Hor{\mathrm{Hor}}
\def\Ad{\mathrm{Ad}}
\def\ad{\mathrm{ad}}
\def\id{\mathrm{id}}

\def\Ker{\mathrm{Ker}}
\def\Im{\mathrm{Im}}
\def\inv{\mathrm{inv}}

\def\c{\mathrm{c}}

\def\Mor{\textsc{Mor}}
\def\N{\mathbb{N}}

\def\triv{\mathrm{triv}}
\def\qtrs{\mathrm{qtrs}}
\def\l{\mathrm{L}}
\def\c{\mathrm{c}}
\def\End{\mathrm{End}}
\def\r{\mathrm{R}}

\def\G{\mathcal{G}}

\def\T{\mathcal{T}}

\newtheorem{Lemma}{Lemma}[section]
\newtheorem{Remark}[Lemma]{Remark}
\newtheorem{Example}[Lemma]{Example}
\newtheorem{Theorem}[Lemma]{Theorem}

\newtheorem{Proposition}[Lemma]{Proposition}

\begin{document}
\date{\today}
\title{Some Constructions on Quantum Principal Bundles}
\author{Gustavo Amilcar Salda\~na Moncada}
\address{Gustavo Amilcar Salda\~na Moncada\\
Mathematics Research Center, CIMAT}
\email{gustavo.saldana@cimat.mx}
\begin{abstract}
This paper works as an appendix of the papers \emph{Geometry of Associated Quantum Vector Bundles and the Quantum Gauge Group}, \emph{Quantum Principal Bundles and Yang--Mills Scalar Matter Fields} and \emph{Yang--Mills--Connes Theory and Quantum Principal $SU(N)$--Bundles over Toric Non--Commutative Manifolds}. Here, we are going to prove the following six statements in the theory of quantum principal bundles:\\
\begin{enumerate}
    \item The universal differential envelope $\ast$--calculus of a matrix (compact) Lie group, for the classical bicovariant First Order Differential $\ast$--Calculus, is the algebra of differential forms.
     \item An example of a quantum principal bundle in which the space of base forms is not generated by the base space.
     \item The group isomorphism between convolution--invertible maps and covariant left module isomorphisms at the level of differential calculus
     \item The way the maps $\{T^V_k \}$ from Remark \ref{rema} look in differential geometry.
     \item The Gauge Principle in non--commutative geometry.
     \item Formal adjointability of finitely generated modules.
     \\
\end{enumerate}
Furthermore, we are going to present an example of the theory of the reference \cite{sald} in a special {\it classical/quantum hybrid} principal bundle.

\end{abstract}
\maketitle

\section{Universal Differential Envelope $\ast$--calculus in Differential Geometry}

A compact matrix quantum group (or a just a quantum group) will be denoted by $\G$, and its dense $\ast$--Hopf (sub)algebra will be denoted by
\begin{equation}
    \label{2.f0}
    H^\infty:=(H,\cdot,\mathbbm{1},\Delta,\epsilon,S,\ast),
\end{equation}
where  $\Delta$ is the coproduct, $\epsilon$ is the counity and $S$ is the coinverse. The space $H^\infty$ shall be treated as the algebra of all {\it polynomial functions} defined on $\G$.  In the same way a (smooth right) $\G$--corepresentation on a $\C$--vector space $V$ is a linear map $$\delta^V: V \longrightarrow V \otimes H$$ such that
\begin{equation} 
\label{2.f1}
(\id_V\otimes \epsilon )\circ \delta^V\cong \id_V \end{equation}
and 
\begin{equation} 
\label{2.f2}
(\id_V\otimes \Delta )\circ \delta^V=(\delta^V\otimes \id_H ) \circ \delta^V.
\end{equation}

\noindent We say that the corepresentation is {\it finite--dimensional} if $\dim_{\C}(V)< |\N|$. $\delta^V$ usually receives the name of the {\it (right) coaction} of $\G$ on $V$. It is worth mentioning that in the general theory (\cite{woro1}), equation (\ref{2.f1}) is not necessary. 

The set of all corepresentation morphisms between two corepresentations $\delta^V$, $\delta^W$ will be denoted as 
\begin{equation}
\label{2.f4}
\Mor(\delta^V,\delta^W),
\end{equation}

A $\G$--corepresentation $\delta^V$ is {\it reducible} if there exists a non--trivial subspace $L$ $(L\not=\{0\}, V)$ such that $\delta^V(L)=L\otimes H$ and $\delta^V$ is {\it unitary} if viewed as an element of $B(V)\otimes H$ (with $B(V):=\{f:V\longrightarrow V\mid f \mbox{ is linear} \}$) is unitary. Of course, for the last definition, it is necessary an inner product $\langle-| -\rangle$ on $V$. In \cite{woro1}, Woronowicz proved that for every finite--dimensional $\G$--corepresentation on $V$ there is an inner product $\langle-|-\rangle$ (not necessarily unique) making the corepresentation unitary. Thus, henceforth we are going to consider that every finite--dimensional $\G$--corepresentation is unitary. In \cite{woro1} one can find a proof of the following theorem.

\begin{Theorem}
\label{rep}
Let $\T$ be a complete set of mutually non--equivalent irreducible unitary (necessarily finite--dimensional) $\G$--corepresentations with $\delta^\C_\triv$ $\in$ $\T$ (the trivial corepresentation on $\C$). For any $\delta^V$ $\in$ $\T$ that coacts on $(V,\langle-|-\rangle)$,
\begin{equation}
\label{2.f8}
\delta^V(e_i)=\sum^{n_{V}}_{j=1} e_j\otimes g^{V}_{ji},
\end{equation}
where $\{ e_i\}^{n_{V}}_{i=1}$ is an ortonormal basis of $V$ and $\{g^{V}_{ij}\}^{n_{V}}_{i,j=1}$ $\subseteq$ $H$. Then $\{g^{V}_{ij}\}_{\delta^V,i,j}$ is a linear basis of $H$, where the index $\delta^V$ runs on $\T$ and $i$, $j$ run from $1$ to $n_{V}=\mathrm{dim}_\C(V)$.
\end{Theorem}

For every $\delta^V$ $\in$ $\T$, the set $\{g^{V}_{ij}\}^{n_{V}}_{ij=1}$ satisfies
\begin{equation}
\label{2.f9}
\phi(g^{\alpha}_{ij})=\sum^{n_{\alpha}}_{k=1}g^{\alpha}_{ik}\otimes g^{\alpha}_{kj},\quad \kappa(g^{\alpha}_{ij})=g^{\alpha\,\ast}_{ji},\quad \sum^{n_{\alpha}}_{k=1}g^{\alpha}_{ik}g^{\alpha\,\ast}_{jk}=\sum^{n_{\alpha}}_{k=1}g^{\alpha\,\ast}_{ki}g^{\alpha}_{kj}=\delta_{ij}\mathbbm{1}, \quad \epsilon(g^{\alpha}_{ij})=\delta_{ij}
\end{equation}
with $\delta_{ij}$ being the Kronecker delta, among other properties \cite{woro1}.

Let $(\Gamma,d)$ be a bicovariant First Order Differential $\ast$--Calculus ($\ast$--FODC) over $\G$ (\cite{woro2,stheve}) and consider the $\C$--vector space given by
\begin{equation}
\label{2.f13}
\mathfrak{qg}^\#:=\{\theta \in \Gamma \mid \Phi_\Gamma(\theta)=\mathbbm{1}\otimes \theta  \}   =\mathbbm{1}\otimes {\Ker(\epsilon)\over \mathcal{R}}\cong {\Ker(\epsilon)\over \mathcal{R}}.
\end{equation}
This space allows to consider the quantum germs map 
\begin{equation}
\label{2.f14}
\begin{aligned}
\pi:H &\longrightarrow \mathfrak{qg}^\#\\
g &\longmapsto [g-\epsilon(g)\mathbbm{1}]_\mathcal{R},
\end{aligned}
\end{equation}
where $[h]_\mathcal{R}$ denotes the equivalence class of the element $h$ $\in$ $H$. The map $\pi$ has several useful properties, for example the restriction map $\pi|_{\Ker(\epsilon)}$ is surjective and 
\begin{equation}
\label{properties}
    \begin{aligned}
\ker(\pi)=\mathcal{R}\oplus \C\mathbbm{1}, & \qquad\;\;\, \pi(g)=S(g^{1})dg^{2},&\quad  dg=g^{(1)}\pi(g^{(2)}),\\  \pi(g)=-(dS(g^{(1)}))g^{(2)},& \quad  dS(g)=-\pi(g^{(1)})S(g^{(2)}),& \quad  \pi(g)^\ast=-\pi(S(g)^\ast)
    \end{aligned}
\end{equation}
for all $g$ $\in$ $H$ \cite{stheve}. It is worth mentioning that the canonical right action ${_\Gamma}\Phi$ of the $\ast$--FODC leaves $\mathfrak{qg}^\#$ invariant (\cite{stheve}); so 
\begin{equation}
\label{2.f15}
\ad:={_\Gamma}\Phi|_{\mathfrak{qg}^\#}: \mathfrak{qg}^\#\longrightarrow \mathfrak{qg}^\#\otimes H
\end{equation}
is a $\G$--corepresentation and it fulfills \cite{stheve}
\begin{equation}
\label{2.f15.1}
\ad\circ \pi= (\pi\otimes \id_H)\circ \Ad,
\end{equation}
where 
\begin{equation}
    \label{2.f9.5}
    \Ad: H\longrightarrow H\otimes H, \qquad g\longmapsto g^{(2)}\otimes S(g^{(1)})\otimes g^{(3)}
\end{equation}
is the (right) adjoint coaction of $H$.

On the other hand there is a right $H$--module structure on $\mathfrak{qg}^\#$ given by 
\begin{equation}
\label{2.f16}
\theta \circ g=S(g^{(1)})\theta g^{(2)}=\pi(hg-\epsilon(h)g)
\end{equation}
for every $\theta=\pi(h)$ $\in$ $\mathfrak{qg}^\#$; and it satisfies \cite{stheve} $$(\theta \circ g)^\ast=\theta^\ast\circ S(g)^\ast.$$

Let $(\Gamma,d)$ be a $\ast$--FODC over $\G$. Consider the space
  \begin{equation}
      \label{udtensor}
      \Gamma^\wedge:=\otimes^\bullet_H\Gamma/\mathcal{Q}\,,\quad  \otimes^\bullet_H\Gamma:=\bigoplus_k (\otimes^k_H\Gamma)\quad \mbox{ with }\quad\otimes^k_H\Gamma:=\underbrace{\Gamma\otimes_H\cdots\otimes_H \Gamma}_{k\; times},
  \end{equation}
  where $\mathcal{Q}$ is the bilateral ideal of $\otimes^\bullet_H\Gamma$ generated by 
  \begin{equation}
      \label{udtensor1}
     \sum_i dg_i\otimes_H dh_i \quad \mbox{ such that } \quad \sum_i g_i\,dh_i=0,
  \end{equation}
  for $g_i$, $h_i$ $\in$ $H$. On the other hand, for a given $t=\vartheta_1\cdots \vartheta_n$ $\in$ $\Gamma^{\wedge\,n }$ the linear map
  \begin{equation}
      \label{diff12}
     d:\Gamma^\wedge\longrightarrow \Gamma^\wedge 
  \end{equation}
  given by $d(t)=d(\vartheta_1\cdots \vartheta_n)=\displaystyle \sum^n_{j=1}(-1)^{j-1}\vartheta_1\cdots \vartheta_{j-1}\cdot d\vartheta_j\cdot \vartheta_{j+1}\cdots \vartheta_n$ $\in$ $\Gamma^{\wedge\,n+1 }$ is well--defined, satisfies the graded Leibniz rule and $d^2=0$ \cite{micho1,stheve}. In other words, $d$ is {\it natural} extension in $\Gamma^\wedge$ of the differential of $\Gamma$. In this way $(\Gamma^\wedge,d, \ast)$ is a graded differential $\ast$--algebra generated by $\Gamma^{\wedge\,0}= H$ (the degree $0$ elements) \cite{micho1,stheve}.

Let $(\Gamma,d)$ be a bicovariant $\ast$--FODC.  Then, the coproduct can be extended to a graded differential $\ast$--algebra morphism
 \begin{equation}
\label{2.f3.5}
\Delta: \Gamma^\wedge \longrightarrow \Gamma^\wedge \otimes  \Gamma^\wedge.
\end{equation}
In particular, in accordance with \cite{micho1,stheve}, we have 
\begin{equation}
    \label{coproduc.1}
    \Delta(\theta)=\mathbbm{1}\otimes \theta+\ad(\theta).
\end{equation}
for all $\theta$ $\in$ $\mathfrak{qg}^\#$.

 The counit and the coinverse can also be extended. In fact,  consider the linear map 
\begin{equation}
\label{2.f3.6}
\epsilon:\Gamma^\wedge\longrightarrow \C 
\end{equation}
defined by $\epsilon|_{H}:=\epsilon$ and $\epsilon|_{\Gamma^{\wedge k}}:=0$ for $k\geq 1$. 

On the other hand, for any $g$ $\in$ $\Ker(\epsilon)$, define $$S^1(g):=-\pi(g^{(2)})S(g^{(3)})S(S(g^{(1)})).$$ Since $\Ad(\mathcal{R})\subseteq \mathcal{R}\otimes H$, we obtain $S^1(g)=0$ for all $g$ $\in$ $\mathcal{R}$. Hence, there exists a well--defined linear map  $S^1:\mathfrak{qg}^\# \longrightarrow \Gamma.$  In accordance with \cite{micho1}, we can extend $S^1$ to the whole $\Gamma$ in such a way that the following equalities hold
 $$S^1(h\,\pi(g))=S^1(\pi(g))\,S(h)\quad \mbox{ and }\quad S^{1}(h\,dg)=d(S(g))\,S(h) $$ for all $g$, $h$ $\in$ $H$. According to \cite{micho1,stheve}, we  extend $S^0:=S$ and $S^1$ to a graded antimultiplicative linear map 
\begin{equation}
\label{2.f3.7}
S:\Gamma^\wedge\longrightarrow \Gamma^\wedge
\end{equation}
which commutes with the differential (\cite{micho1}).  These maps define a graded differential $\ast$--Hopf algebra structure 
\begin{equation}
\label{2.f3.8}
\Gamma^{\wedge\,\infty}:=(\Gamma^\wedge,\cdot,\mathbbm{1},\Delta,\epsilon, S,d,\ast),
\end{equation}
on $\Gamma^\wedge$ which extends $H^\infty=(H,\cdot,\mathbbm{1},\Delta,\epsilon, S,d,\ast)$ \cite{micho1}.  

Now it is possible to consider the right adjoint coaction of $\Gamma^\wedge$ by taking
\begin{equation}
\label{2.f11}
\Ad:\Gamma^\wedge \longrightarrow \Gamma^\wedge \otimes \Gamma^\wedge
\end{equation}
such that $$\Ad(t)=(-1)^{\partial t^{(1)}\partial t^{(2)}} t^{(2)}\otimes S(t^{(1)})t^{(3)},$$ where $\partial x$ denotes the grade of $x$ and   $(\id_{\Gamma^\wedge}\otimes \Delta)\Delta(t)=(\Delta\otimes \id_{\Gamma^\wedge})\Delta(t)=t^{(1)}\otimes t^{(2)}\otimes t^{(3)}.$  Clearly, $\Ad$ extends the right adjoint coaction $\Ad$ of $H$.

Let us define
\begin{equation}
    \label{2.f12.1}
    \mathfrak{qg}^{\#\wedge}=\otimes^\bullet \mathfrak{qg}^\#/A^\wedge, \quad \otimes^\bullet\mathfrak{qg}^\#:=\oplus_k (\otimes^k\mathfrak{qg}^\#)\; \mbox{ with }\;\otimes^k\mathfrak{qg}^\#:=\underbrace{\mathfrak{qg}^\#\otimes\cdots\otimes \mathfrak{qg}^\#}_{k\; times},
\end{equation}
where $A^\wedge$ is the bilateral ideal of $\otimes^\bullet\mathfrak{qg}^\#$ generated by $$\pi(g^{(1)})\otimes \pi(g^{(2)})$$  with $g$ $\in$ $\mathcal{R}$ \cite{micho1}. The space $\mathfrak{qg}^{\#\wedge}$ is a graded differential $\ast$--subalgebra (\cite{micho1}) and it satisfies $\mathfrak{qg}^{\#\wedge}=\{ \theta \in \Gamma^{\wedge}\mid \Phi_{\Gamma^{\wedge}}(\theta)=\mathbbm{1}\otimes \theta\}$. Also one can extend the right $H$--module structure of $\mathfrak{qg}^\#$ (see equation (\ref{2.f16})) to $\mathfrak{qg}^{\#^\wedge}$ by means of 
\begin{equation}
    \label{2.f12.4}
    1\circ g=\epsilon(g),\quad (\theta_1\theta_2)\circ g=(\vartheta_1\circ g^{(1)})(\theta_2\circ g^{(2)}).
\end{equation}
It is worth mentioning that 
\begin{equation}
    \label{2.f12.2}
    d\pi(g)=-\pi(g^{(1)})\pi(g^{(2)}).
\end{equation}
According to \cite{micho1} the following identification holds:
\begin{equation}
    \label{2.f12.3}
    \Gamma^{\wedge}= H\otimes \mathfrak{qg}^{\#\wedge}.
\end{equation}

The following example is one of the purposes of this paper
 
\begin{Example}
    \label{e.0}
Let $G$ be a compact matrix Lie group ($G\subseteq M_k(\C))$), and let $\G$ be its associated quantum group \cite{woro1}. If 
    \begin{equation}
        \label{2.2.f1}
        H^\infty=(H,\cdot,\mathbbm{1},\Delta,\epsilon,S,\ast,)
    \end{equation}
     is its dense $\ast$--Hopf algebra, then $H$ is the $\ast$--algebra generated by the smooth $\C$--valued functions
    \begin{equation}
    \label{2.2.f2}
    \begin{aligned}
    w_{ij}: G &\longrightarrow \C \\   
          A &\longmapsto a_{ij}.
    \end{aligned}
\end{equation}
    with $A=(a_{ij})$. The coproduct of an element $g$ $\in$ $H$ is defined by 
    \begin{equation}
        \label{classicalcoproduct}
        \Delta(g): G\times G \longmapsto \C, \qquad \Delta(g)(A,B)=g(AB)\quad \mbox{ with }\quad A,\,B \,\in \, G.
    \end{equation}
In particular, we have $\Delta(w_{ij})= \displaystyle \sum^n_{k=1} w_{ik}\otimes w_{kj}.$ The counit and the coinverse are defined as follows: 
\begin{equation}
    \label{classicalcounit}
    \epsilon:H \longrightarrow \C,\qquad g \longmapsto g(e),
\end{equation}
where $e$ $\in$ $H$ is the identity element; and 
\begin{equation}
    \label{classicalantipode}
    S(g): G \longrightarrow \C,\qquad A \longmapsto S(g)(A)=g(A^{-1})
\end{equation}
for all $g$ $\in$ $H$. 

    Since $G$ is parallelizable, we have the identification $G\times \mathfrak{g}=TG=\mathfrak{g}\times G$, where $\mathfrak{g}$ is the Lie algebra of $G$. This diffeomorphism can be given by 
    \begin{equation}
        \label{prod1}
        (A,v) \longmapsto (d'L_A)_e(v),
    \end{equation}
    where $(d'L_A)_e$ is the usual de Rham differential $d'$ at $e$ of the left product on $G$ with $A$; and by 
   \begin{equation}
        \label{prod2}
        (v,A) \longmapsto (d'R_A)_e(v),
    \end{equation}
    where $(d'R_A)_e$ is the de Rham differential $d'$ at  $e$ of the right product on $G$ with $A$.  

   Consider $d=d'$. Under equation (\ref{prod1}) we have $dg=g^{(1)}\otimes dg^{(2)}_e$ for all $g$ $\in$ $H$. By defining $$\mathfrak{h}:=\mathrm{span}_\C\{dg_e\mid g\, \in\, H \} $$ we obtain that 
    $$(\Gamma=H\otimes \mathfrak{h},d) $$ is a $\ast$--FODC over $\G$. Moreover, it is left--covariant under the map $\Phi_\Gamma=\Delta\otimes \id_{\mathfrak{h}}$. 

    On the other hand, under equation (\ref{prod2}), we have $dg=dg^{(1)}_e\otimes g^{(2)}$ for all $g$ $\in$ $H$. Thus $$(\Gamma,d)\cong(\mathfrak{h}\otimes H,d) $$ and it follows that $(\Gamma,d)$ is also right--covariant by the map ${_\Gamma}\Phi=\id_{\mathfrak{h}}\otimes\Delta.$ Hence, $(\Gamma,d)$ is bicovariant. According to  \cite{micho1,stheve}, there exists a right $H$--ideal $\mathcal{R} \subseteq \Ker(\epsilon)$ such that
    $$\{\theta \in \Gamma \mid \Phi_\Gamma(\theta)=\mathbbm{1}\otimes \theta  \}=\mathbbm{1}\otimes {\Ker(\epsilon)\over \mathcal{R}}= {\Ker(\epsilon)\over \mathcal{R}}=:\mathfrak{qh}^\#$$ However, since $\mathfrak{h}\cong \mathbbm{1}\otimes \mathfrak{h}= \{\theta \in \Gamma \mid \Phi_\Gamma(\theta)=\mathbbm{1}\otimes \theta  \}$ it follows that $$\mathfrak{h}=\mathfrak{qh}^\#$$ We claim that $\mathcal{R}=\Ker^2(\epsilon):=\{\displaystyle \sum^n_{i=1}a_i\,b_i \mid a_i,\, b_i \,\in\, \Ker(\epsilon) \, \mbox{ for some }\, n\,\in\, \N  \}$. In fact, consider the (restricted) quantum germs map
\begin{equation*}
    \begin{aligned}
\pi: \Ker(\epsilon) &\longrightarrow {\Ker(\epsilon)\over \mathcal{R}} \\   
          g &\longmapsto [g]_{\mathcal{R}}.
    \end{aligned}
\end{equation*}
Then, by equation (\ref{properties}), we obtain  $\pi(g)=S(g^{(1)})dg^{(2)}=dg_e$. If $\mathfrak{g}^\#_\C$ is the complexification of the dual space of the Lie algebra $\mathfrak{g}$ of $G$, 
it is well--known that for $X=\{f\in C^\infty_\C(G)\mid \epsilon(f)=0 \}$ we have 
$$\mathfrak{g}^\#_\C={X \over X^2 } $$ with $X^2=\{\displaystyle \sum^n_{i=1}a_i\,b_i \mid a_i,\, b_i \,\in\, X \, \mbox{ for some }\, n\,\in\, \N  \}$ \cite{war}.  Moreover, for every $f$ $\in$ $X$, $$df_e=[f]_{X^2},$$ where $[f]_{X^2}$ denotes the equivalence class of $f$ in $X/X^2$ \cite{war}. In this way,  for every $g$ $\in$ $\Ker(\epsilon)$ we get $dg_e=[g]_\mathcal{R}=[g]_{X^2}$. By this relation,  $g$ $\in$ $\mathcal{R}$ if and only if $g$ $\in$ $X^2\cap \Ker(\epsilon)=\Ker^2(\epsilon),$ which proves that $\mathcal{R}=\Ker^2(\epsilon)$.

By definition, $\mathfrak{h}\subseteq \mathfrak{g}^\#_\C$. Since $\mathfrak{g}_\C \subseteq M_k(\C)$, there exists a linear basis $\{ v_i \}$ of $\mathfrak{h}_\C$  composed of linear combinations of the linear basis $\{E_{lj} \}$ of $M_k(\C)$ ($E_{lj}$ is the matrix with $1$ in the $(l,j)$ position and $0$ elsewhere). Then $(dw_{lj})_e(v_s)$ corresponds to the $(l,j)$ position of the matrix $v_s$ . Consequently, there exist elements $\{(dg_i)_e\}$, which are linear combinations of $(dw_{ij})_e$ $\in$ $\mathfrak{g}$, such that $(dg_i)_e(v_j)=\delta_{ij}$ with $\delta_{ij}$ being the Kronecker delta. Thus, $\mathfrak{g}^\#_\C \subseteq \mathfrak{h}$ and hence $\mathfrak{g}^\#_\C = \mathfrak{h}$. This implies that  
\begin{equation}
    \label{2.2.f3}
    (\Gamma,d)=(H\otimes \mathfrak{g}^\#_\C,d);
\end{equation}
so $(\Gamma,d)$ is a  $\ast$--subFODC of the $\ast$--FODC $$(C^\infty_\C(G)\otimes \mathfrak{g}^\#_\C,d=d')$$ over $C^\infty_\C(G)$  of $\C$--valued differential $1$--forms of $G$. Notice that it suffices to take convergent sequences $\{ g_i\}^\infty_{i=1}$ $\subseteq$ $H$ in $(\Gamma,d)$ to recover $(C^\infty_\C(G)\otimes \mathfrak{g}^\#_\C,d)$. It is worth mentioning that for the specific form  $\Gamma=H\otimes \mathfrak{g}^\#_\C$, the right covariance map ${_\Gamma}\Phi$ is given by \cite{stheve} $${_\Gamma}\Phi(g\otimes \pi(h))=g^{(1)}\otimes \pi(h^{(2)})\otimes g^{(2)}S(h^{(1)})h^{(1)}$$ for all $g$, $h$ $\in$ $H$.

 By the Leibniz rule, we have $$d(hg)_e=\pi(hg)=\pi(h)\epsilon(g)+\epsilon(h)\pi(g)$$ for all $h$, $g$ $\in$ $H$ and the right $H$--module structure on $\mathfrak{g}^\#_\C$,  given by  $$\pi(h)\circ g=\pi(hg-\epsilon(h)g),$$ simplifies to 
 \begin{equation}
     \label{2.2.f4}
     \pi(h)\circ g=\epsilon(g)\pi(h).
 \end{equation}
 
Now, consider the universal differential envelope $\ast$--calculus $$(\Gamma^\wedge,d,\ast)$$ of $(\Gamma,d)$. If $g=a\,b$ $\in$ $\Ker^2(\epsilon)$, we obtain
\begin{eqnarray}
\pi(g^{(1)})\otimes \pi(g^{(2)})=\pi(a^{(1)}b^{(1)})\otimes \pi(a^{(2)}b^{(2)})&= &\pi(a^{(1)})\circ b^{(1)}\otimes \pi(a^{(2)})\circ b^{(2)} \nonumber
  \\
&+&
\epsilon(a^{(1)})\pi(b^{(1)})\otimes \pi(a^{(2)})\circ b^{(2)}\nonumber
  \\
  &+ &
  \pi(a^{(1)})\circ b^{(1)}\otimes \epsilon(a^{(2)})\pi(b^{(2)}) \label{2.2.cuentas}
  \\
&+&
\epsilon(a^{(1)})\pi(b^{(1)})\otimes\epsilon(a^{(2)})\pi(b^{(2)})\nonumber
  \\
  &= &
  \pi(b)\otimes \pi(a)+\pi(a)\otimes \pi(b).\nonumber
\end{eqnarray}
and therefore (see equation (\ref{2.f12.1}))
$$\mathfrak{g}^{\#\,\wedge}_\C=\otimes^\bullet \mathfrak{g}^\#_\C/A^\wedge=\bigwedge \mathfrak{g}^\#_\C,$$ where $\displaystyle \bigwedge \mathfrak{g}^\#_\C$ denotes the exterior algebra of $\mathfrak{g}^\#_\C$. Thus (see equation (\ref{2.f12.3}))
\begin{equation}
    \label{2.2.f5}
    (\Gamma^\wedge=H\otimes \bigwedge \mathfrak{g}^\#_\C,d,\ast).
\end{equation}
By the axioms of the exterior derivative, we can conclude that $d: \Gamma^\wedge \longrightarrow \Gamma^\wedge$ is precisely the exterior derivative. Hence $(\Gamma^\wedge,d,\ast)$ is a graded differential $\ast$--subalgebra of the graded differential $\ast$--algebra 
\begin{equation}
    \label{2.2.f5.1}
    (\Omega^\bullet_\C(G)= C^\infty_\C(G)\otimes \bigwedge \mathfrak{g}^\#_\C,d,\ast)
\end{equation}
of $\C$--valued differential forms of $G$. As before, it suffices to consider convergent sequences $\{ g_i\}^\infty_{i=1}$ $\subseteq$ $H$ in $(\Gamma^\wedge,d,\ast)$ to obtain $(\Omega^\bullet_\C(G),d,\ast).$ 

Now we are going to proceed characterizing the graded differential $\ast$--Hopf algebra structure $$ \Gamma^{\wedge\,\infty}=(\Gamma^{\wedge},\cdot,\mathbbm{1},\Delta,\epsilon,S,d,\ast).$$

Let $\theta$ $\in$ $\mathfrak{g}^\#_\C$. Then $\theta=\pi(g)=(dg)_e$ for some $g$ $\in$ $H$. Taking $h$ $\in$ $H$, by equations (\ref{2.f15.1}), (\ref{coproduc.1}) and the fact that $\Delta$ is multiplicative, we have 
\begin{eqnarray*}
    \Delta(h\,\pi(g))&=&h^{(1)}\otimes h^{(2)}\,\pi(g)+(h^{(1)}\otimes h^{(2)})\,\ad(\pi(g))
    \\
    &=& 
    h^{(1)}\otimes h^{(2)}\,\pi(g)+h^{(1)}\,\pi(g^{(2)})\otimes h^{(2)}\,S(g^{(1)})g^{(3)} \;\in\; H\otimes \Gamma\oplus \Gamma\otimes H. 
\end{eqnarray*}
Let $(A,(B,v))$ $\in$ $G\times TG_\C$, $((D,u),C)$ $\in$ $TG_\C\times G$ with $TG_\C=G\times \mathfrak{g}_\C$, for some $v$, $u$ $\in$ $\mathfrak{g}_\C$. Thus
\begin{eqnarray*}
    \Delta(h\,\pi(g))((A,(B,v))\oplus ((D,u),C)&=&
    (h^{(1)}\otimes h^{(2)}\,\pi(g))(A,(B,v))
    \\
    &+&
    (h^{(1)}\,\pi(g^{(2)})\otimes h^{(2)}\,S(g^{(1)})g^{(3)})((D,u),C)
    \\
    &=&
    h^{(1)}(A)\,h^{(2)}(B)\,\pi(g)(v)
    \\
    &+&
    h^{(1)}(D)\,\pi(g^{(2)})(u)\,h^{(2)}(C)\,g^{(1)}(C^{-1})\,g^{(3)}(C)
    \\
    &=&
    h(AB)\,\pi(g)(v)
    \\
    &+&
    h(DC)\,g^{(1)}(C^{-1})\,\pi(g^{(2)})(u)\,g^{(3)}(C)
    \\
    &=&
    h(AB)\,\pi(g)(v)
    \\
    &+&
    h(DC)\,g^{(1)}(C^{-1})\,(dg^{(2)})_e(u)\,g^{(3)}(C) 
    \\
    &=&
    h(AB)\,\pi(g)(v)
    \\
    &+&
    h(DC)\,\left.{d \,g(C^{-1}\,\mathrm{exp}(tu)\,C)\over dt} \right|_{t=0}.
\end{eqnarray*}

According to reference \cite{nodg}, the group product $$m_{TG}:=dm_G:TG_\C\times TG_\C\longrightarrow TG_\C $$  of $TG_\C=G\times \mathfrak{g}_\C$ is given by  $$m_{TG}((R,x),(S,y))=(RS,y+\ad^{\mathrm{class}}_S(x)),$$ where $\ad^{\mathrm{class}}_S$ is the differential at $e$ of the {\it classical} right adjoint action $\Ad^{\mathrm{class}}_S$ by $S$ $\in$ $G$ on $G$, which is given by $$\Ad^{\mathrm{class}}:G\times G\longrightarrow G,\qquad \Ad^{\mathrm{class}}(R,S):=\Ad^{\mathrm{class}}_S(R):=S^{-1}RS;$$ so 
\begin{equation}
    \label{adjointclassical}
    \ad^{\mathrm{class}}:\mathfrak{g}_\C\times G\longrightarrow \mathfrak{g}_\C,\qquad \ad^{\mathrm{class}}(v,S):=\ad^{\mathrm{class}}_S(v)=(d\,\Ad^{\mathrm{class}}_S)_e(v).
\end{equation}
Extending $m_{TG}$ to $$m_{TG}: G\times TG_\C\longrightarrow TG_\C,\qquad \mbox{ and }  \qquad  m_{TG}: TG_\C \times G\longrightarrow TG_\C $$ we get $$m_{TG}(A,(B,v))=(AB,v),\qquad m_{TH}((D,u),C)=(DC,\ad^{\mathrm{class}}_C(u)).$$ Thus
\begin{eqnarray*}
    ((h\,\pi(g))\circ m_{TG})((A,(B,v))\oplus ((D,u),C))&=&(h\,\pi(g))(AB,v)
    \\
    &+&
    (h\,\pi(g))(DC,\ad^{\mathrm{class}}_C(u))
    \\
    &=&
    h(AB)\,\pi(g)(v)
    \\
    &+&
    h(DC)\pi(g)(\ad^{\mathrm{class}}_C(u))
    \\
    &=&
    h(AB)\,\pi(g)(v)
    \\
    &+&
    h(DC)\,\left.{d \,g(C^{-1}\,\mathrm{exp}(tu)\,C)\over dt} \right|_{t=0}.
\end{eqnarray*}
Since every element $\vartheta$ of $\Gamma$ is the sum of elements of the form $h\,\pi(g)$, we conclude that $$\Delta(\vartheta)=\vartheta\circ m_{TG}.$$ By equation (\ref{classicalcoproduct}), the last result and the facts that $(\Gamma^\wedge,d,\ast)$ is generated by its degree $0$ elements and $\Delta$ is a graded differential $\ast$--algebra morphism, it immediately follows that, for all $\vartheta$ $\in$ $\Gamma^\wedge$
\begin{equation}
    \label{classicalcoproduct1}
    \Delta(\vartheta)=\vartheta\circ m_{TG}.
\end{equation}
In light of last equation, we can extend the coproduct $\Delta:\Gamma^\wedge\longrightarrow \Gamma^\wedge\otimes \Gamma^\wedge$ to 
\begin{equation}
    \label{classicalcoproduct2}
    \Delta: \Omega^\bullet_\C(G)\longrightarrow \Omega^\bullet_\C(G)\,\widetilde{\otimes}\,\Omega^\bullet_\C(G),\qquad \vartheta\longmapsto \vartheta\circ m_{TG},
\end{equation}
where $\widetilde{\otimes}$ denotes the corresponding completion of the algebraic tensor product $\otimes$ of graded differential $\ast$--algebras.

On the other hand, let $\pi(g)$ $\in$ $\mathfrak{g}^\#_\C$   and $h$ $\in$ $H$. In light of the definition of the antipode $S$ in equations (\ref{2.f3.7}), (\ref{classicalantipode}), we have 
\begin{eqnarray*}
    S(h\,\pi(g))=S(\pi(g))\,S(h)&=&-\pi(g^{(2)})\,S(g^{3})\,S(S(g^{(1)}))\,S(h)
    \\
    &=&
    -\pi(g^{(2)})\,S(g^{3})\, g^{(1)}\,S(h)
    \\
    &=&
    -S(g^{3})\, g^{(1)}\,S(h)\,\pi(g^{(2)})\;\in\; \Gamma
\end{eqnarray*}
Let $(A,v)$ $\in$ $TG_\C=G\times \mathfrak{g}_\C$. Then
\begin{eqnarray*}
    S(h\,\pi(g))(A,v)&=&-(S(g^{3})\, g^{(1)}\,S(h)\,\pi(g^{(2)}))(A,v)
    \\
    &=&g^{(3)}(A^{-1})\,g^{(1)}(A)\,h(A^{-1})\,\pi(g^{(2)})(v)
    \\
    &=&
    -h(A^{-1})\,g^{(1)}(A)\,\pi(g^{(2)})(v)\,g^{(3)}(A^{-1})
    \\
    &=&
    -h(A^{-1})\,g^{(1)}(A)\,(dg^{(2)})_e(v)\,g^{(3)}(A^{-1})
    \\
    &=&
    -h(A^{-1})\,\left.{d \,g(A\,\mathrm{exp}(tv)\,A^{-1})\over dt} \right|_{t=0}.
\end{eqnarray*}

\noindent According to \cite{nodg}, the map $$\iota_{TG}:=d\iota_G: TG_\C\longrightarrow TG_\C$$ that sends an element of $TG_\C=G\times \mathfrak{g}_\C$ to its multiplicative inverse,  is given by
$$\iota_{TG}(A,v)=(A,v)^{-1}=(A^{-1},-\ad^{\mathrm{class}}_{A^{-1}}(v)).$$ In this way, we have 
\begin{eqnarray*}
    ((h\,\pi(g))\circ \iota_{TG})(A,v)&=&(h\,\pi(g))(A^{-1},-\ad^{\mathrm{class}}_{A^{-1}}(v))
    \\
    &=&
    -h(A^{-1})\,\pi(g)(\ad^{\mathrm{class}}_{A^{-1}}(v))
    \\
    &=&
    -h(A^{-1})\,\left.{d \,g(A\,\mathrm{exp}(tv)\,A^{-1})\over dt} \right|_{t=0}. 
\end{eqnarray*}
Since every element $\vartheta$ of $\Gamma$ is the sum of elements of the form $h\,\pi(g)$, we conclude that $$S(\vartheta)=\vartheta\circ \iota_{TG}.$$ By equation (\ref{classicalantipode}), the last result and the facts that $(\Gamma^\wedge,d,\ast)$ is generated by its degree $0$ elements and $S$ is a graded antimultiplicative linear map that commutes with the differential, it immediately follows that, for all $\vartheta$ $\in$ $\Gamma^\wedge$
\begin{equation}
    \label{classicalantipode1}
    S(\vartheta)=\vartheta\circ \iota_{TG}.
\end{equation}
In light of last equation, we can extend the coantipode $S:\Gamma^\wedge\longrightarrow \Gamma^\wedge$ to
\begin{equation}
    \label{classicalantipode2}
    S:\Omega^\bullet_\C(G)\longrightarrow \Omega^\bullet_\C(G), \qquad \vartheta\longmapsto \vartheta\circ \iota_{TG}.
\end{equation}

Finally, according to \cite{nodg}, the multiplicative neutral element of $TG_\C=G\times \mathfrak{g}_\C$ is $e\cong (e,0)$; thus (see equation (\ref{2.f3.6}))
$$\vartheta(e,0)=0=\epsilon(\vartheta)$$ for all $\vartheta$ $\in$ $\Gamma$. By the last result and the facts that $(\Gamma^\wedge,d,\ast)$ is generated by its degree $0$ elements and $\epsilon$ is a graded differential $\ast$--algebra morphism, it immediately follows that  $$\epsilon: \Gamma^{\wedge}\longrightarrow \C $$ is the evaluation at $e$ $\in$ $G$. As before, we can extend the counit $\epsilon$ to 
\begin{equation}
    \label{classicalcounit1}
    \epsilon:\Omega^\bullet_\C(G)\longrightarrow \Omega^\bullet_\C(G), \qquad \vartheta\longmapsto \vartheta(e).
\end{equation}

Hence, we conclude that the structure of graded differential $\ast$--Hopf algebra of $\Gamma^\wedge$ reflects the fact that for a given Lie group, its complexified tangent bundle $TG_\C=G\times \mathfrak{g}_\C$ is also a Lie group (\cite{nodg}).    

It is worth mentioning that the previous calculations also show that the (right) adjoint $H$--coaction on $H$ (see equation (\ref{2.f9.5})) is the pull--back of the {\it classical} (right) adjoint action of $G$ on $G$, and the (right) adjoint $H$--coaction on $\mathfrak{qh}^\#=\mathfrak{g}^\#_\C$ (see equation (\ref{2.f15})) is the pull--back of the {\it classical} (right) adjoint action of $G$  on $\mathfrak{g}_\C$. In other words
\begin{equation}
    \label{adstructure1}
        \Ad(g)=g\circ \Ad^{\mathrm{class}} 
\end{equation}
for all $g$ $\in$ $H$, and
\begin{equation}
    \label{adstructure1.1}
    \ad(\theta)=\theta\circ \ad^{\mathrm{class}}
\end{equation}
for all $\theta$ $\in$ $\mathfrak{g}^\#_\C$. 

Consider the {\it classical} adjoint Lie algebra representation of  $\mathfrak{g}_\C$
\begin{equation}
    \label{adstructure2}
    \mathrm{c}: \mathfrak{g}_\C\times \mathfrak{g}_\C\longrightarrow \mathfrak{g}_\C
\end{equation}
which is given by 
\begin{eqnarray*}
    \mathrm{c}(v,u)&=&\left.{d\,\ad^{\mathrm{class}}_{\mathrm{exp}(tu)}(v)\over dt }\right|_{t=0}
    \\
    &=&
    \left.{d\,(d\,\Ad^{\mathrm{class}}_{\mathrm{exp}(tu)})_e(v)\over dt }\right|_{t=0}
    \\
    &=&
    \left.{d\over dt }\left(\left.{d\over dx }(\mathrm{exp}(tu)^{-1}\,\mathrm{exp}(xv)\,\mathrm{exp}(tu))\right|_{x=0}\right)\right|_{t=0}
    \\
    &=&
    [v,u]_\C, 
\end{eqnarray*}
where $[-,-]_\C$ is the complexified Lie bracket. 

Now, consider the linear map
\begin{equation}
    \label{adstructure3}
    \mathrm{c}^T:=(\id_{\mathfrak{g}^\#_\C}\otimes \pi)\circ \ad: \mathfrak{g}^\#_\C\longrightarrow \mathfrak{g}^\#_\C\otimes \mathfrak{g}^\#_\C
\end{equation}
given by $$\mathrm{c}^T(\theta)=\mathrm{c}^T(\pi(g))= \pi(g^{(2)})\otimes \pi(S(g^{(1)})g^{(3)})$$ for all $\theta$ $\in$ $\mathfrak{g}^\#_\C$.

\noindent Then, for all $v$, $u$ $\in$ $\mathfrak{g}_\C$, we have

\begin{eqnarray*}
    (\mathrm{c}^T(\theta))(v,u)&=&( \pi(g^{(2)})\otimes \pi(S(g{(1)})g^{(3)}))(v,u)
    \\
    &=&
    \pi(g^{(2)})(v)\,\pi(S(g{(1)})g^{(3)})(u)
    \\
    &=&
   \left( \left.{dg^{(2)}(\mathrm{exp}(xv))\over dx}\right|_{x=0}\right) \left(\left.{d\over dt}(g^{(1)}(\mathrm{exp}(tu)^{-1})\,g^{(3)}(\mathrm{exp}(tu))\right|_{t=0}\right)
    \\
    &=&
    \left.{d\over dt}(g^{(1)}(\mathrm{exp}(tu)^{-1})\,\left.{dg^{(2)}(\mathrm{exp}(xv))\over dx}\right|_{x=0}\,g^{(3)}(\mathrm{exp}(tu))\right|_{t=0}
    \\
    &=&
    \left.{d\over dt}(g^{(1)}(\mathrm{exp}(tu)^{-1})\,(dg^{(2)})_e(v)\,g^{(3)}(\mathrm{exp}(tu))\right|_{t=0}
    \\
    &=&
    \left.{d\over dt}\left( \left.{d\over dx} g(\mathrm{exp}(tu)^{-1}\,\mathrm{exp}(xv)\,\mathrm{exp}(tu) )\right|_{x=0}\right)\right|_{t=0}
    \\
    &=&
     (dg)_e\left( \left.{d\over dt }\left(\left.{d\over dx }(\mathrm{exp}(tu)^{-1}\,\mathrm{exp}(xv)\,\mathrm{exp}(tu))\right|_{x=0}\right)\right|_{t=0}\right)
    \\
    &=&
    (dg_e)(c(v,u))
    \\
    &=&
    \pi(g)(c(v,u))=\theta(c(v,u)).
\end{eqnarray*}
We conclude that
\begin{equation}
    \label{adstructure4}
    \mathrm{c}^T(\theta)=\theta\circ \mathrm{c}.
\end{equation}
for all $\theta$ $\in$ $\mathfrak{g}^\#_\C$.

The example developed  explicitly shows that the universal differential envelope $\ast$--calculus $(\Gamma^\wedge,d,\ast)$ is a proper generalization of the algebra $(\Omega^\bullet_\C(G),d,\ast)$ of $\C$--valued differential forms of $G$ in non--commutative geometry.  In this way, for a given quantum group $\G$ and a bicovariant $\ast$--FODC $(\Gamma,d)$ over a quantum group $\G$, the triplet $$(\Gamma^\wedge,d,\ast)$$ will be interpreted as the $\ast$--algebra of {\it quantum differential forms} of $\G$. In this sense,  the quantum dual Lie algebra 
$$\mathfrak{qh}^\#={\Ker(\epsilon)\over \mathcal{R}}$$ plays the role of the {\it dualization} $\mathfrak{g}^\#_\C$ of $\mathfrak{g}_\C$  and the  $H$--corepresentation $\ad$ plays the role of the {\it dualization} of the right adjoint action $\ad^{\mathrm{class}}$ of $G$ on $\mathfrak{g}_\C$ (see equation (\ref{adstructure1.1})). Furthermore, by equation (\ref{adstructure4}), the linear map
\begin{equation}
    \label{trans1}
    \mathrm{c}^T:= (\id_{\mathfrak{qh}^\#}\otimes \pi)\circ \ad : \mathfrak{qh}^\#\longrightarrow \mathfrak{qh}^\#\otimes \mathfrak{qh}^\# 
\end{equation}
will be interpreted as the {\it quantum Lie bracket}. 
\end{Example}

\section{An example about the space of base forms.}

Let $(B,\cdot,\mathbbm{1},\ast)$ be a quantum space and let $\G$ be a quantum group. A {\it quantum principal $\G$--bundle} over $B$ (abbreviated "qpb") is a quantum structure formally represented by the triplet 
\begin{equation}
\label{2.f17}
\zeta=(P,B,\Delta_P),
\end{equation}
where $(P,\cdot,\mathbbm{1},\ast)$ is called the {\it quantum total space}, and $(B,\cdot,\mathbbm{1},\ast)$ a quantum subspace, which receives the name {\it quantum base space}. Furthermore, $$\Delta_P:P \longrightarrow P\otimes H$$ is a $\ast$--algebra morphism that satisfies
\begin{enumerate}
\item $\Delta_P$ is a $\G$--corepresentation.
\item $\Delta_P(x)=x\otimes \mathbbm{1}$ if and only if $x$ $\in$ $B$.
\item The linear map $\beta:P\otimes P\longrightarrow P\otimes H$ given by $$\beta(x\otimes y):=x\cdot \Delta_P(y):=(x\otimes \mathbbm{1})\cdot \Delta_P(y) $$ is surjective. 
\end{enumerate}
 
In general, one does not need a quantum group, a $\ast$--Hopf algebra suffices \cite{micho2,stheve}. Given $\zeta$ a qpb over $B$, a {\it differential calculus} on it is:
 \begin{enumerate}
 \item A graded differential $\ast$--algebra $(\Omega^\bullet(P),d,\ast)$ generated by $\Omega^0(P)=P$ ({\it quantum differential forms of $P$}).
 \item  A bicovariant $\ast$--FODC $(\Gamma,d)$ over $\G$.
 \item The map $\Delta_P$ is extendible to a graded differential $\ast$--algebra morphism $$\Delta_{\Omega^\bullet(P)}:\Omega^\bullet(P)\longrightarrow \Omega^\bullet(P)\otimes \Gamma^{\wedge}.$$ Here we have considered that $\otimes$ is the tensor product of graded differential $\ast$--algebras.
 \end{enumerate}

Note that if $\Delta_{\Omega^\bullet(P)}$ exists, it is unique because our graded differential $\ast$--algebras are generated by its degree $0$ elements. Furthermore,  $\Delta_{\Omega^\bullet(P)}$ is a graded differential $\Gamma^\wedge$--corepresentation on $\Omega^\bullet(P)$ \cite{micho2}. In this setting, the space of horizontal forms is defined as
\begin{equation}
\label{2.f18}
\Hor^\bullet P\,:=\{\varphi \in \Omega^\bullet(P)\mid \Delta_{\Omega^\bullet(P)}(\varphi)\, \in \, \Omega^\bullet(P)\otimes H \},
\end{equation}
and it is a graded  $\ast$--subalgebra of $\Omega^\bullet(P)$ \cite{stheve}. Since $\Delta_{\Omega^\bullet(P)}(\Hor^\bullet P)\subseteq \Hor^\bullet P\otimes H,$ the map 
\begin{equation}
\label{2.f19}
\Delta_\Hor:=\Delta_{\Omega^\bullet(P)}|_{\Hor^\bullet P}: \Hor^\bullet P \longrightarrow \Hor^\bullet P\otimes H
\end{equation}
is $\G$--corepresentation on $\Hor^\bullet P$. Also, one can define the space of {\it base} forms ({\it quantum differential forms of $B$}) as 
\begin{equation}
\label{2.f20}
\Omega^\bullet(B):=\{\mu \in \Omega^\bullet(P)\mid \Delta_{\Omega^\bullet(P)}(\mu)=\mu\otimes \mathbbm{1}\}.
\end{equation}
The space of base forms is a graded differential $\ast$--subalgebra of $(\Omega^\bullet(P),d,\ast)$, and in general, it is not generated by $B$. In fact, the following example, one of the main points of this paper, illustrates a situation in which $(\Omega^\bullet(B),d,\ast)$ is not generated by $B$.
\begin{Example}
    \label{e.2}
    Let $\G$ be the quantum group associated with $U(1)$ and let $$(H=\C[z,z^{\ast}],\cdot,\mathbbm{1},\Delta,\epsilon,S,\ast)$$ be its dense $\ast$--Hopf algebra. Consider $B:=\{ \lambda\,\mathbbm{1}\mid\lambda \in \C \}$ and define the triplet $$\zeta=(P:=H,B,\Delta_P:=\Delta),$$ which is a quantum principal $\G$--bundle over $B$. 

    For quantum differential forms of $P$,
    we will use the universal differential envelope $\ast$--calculus $(\Omega^\bullet(P),d,\ast)$ presented in Example \ref{e.0} for $U(1)$. 
     Recall that $${_\inv}\Omega^1(P):={\Ker(\epsilon)\over \Ker^2(\epsilon) }=\mathfrak{u}(1)^\#_\C$$ is the complexification of the dual space of $\mathfrak{u}(1)$, the Lie algebra of $U(1)$. Of course in this case, $\Omega^\bullet(P)$ has no elements of degree $n \geq 2$, and $${_\inv}\Omega^1(P)=\mathrm{span}_\C\{\pi(z)=dz_e \}.$$  In accordance with \cite{stheve}, the set $\{\pi(z) \}$ is a left $P$--basis of $\Omega^1(P)$. 
    
    For quantum differential forms of $\G$, we will take 
    the universal differential envelope $\ast$--calculus of the bicovariant  $\ast$--FODC over $\G$ defined by  $$ \mathcal{R}=\Ker(\epsilon).$$  In this case, $\Gamma=H\otimes \{0\}=\{ 0\}$, so $\Gamma^\wedge=\Gamma^{\wedge\,0}=H$, and $d=0$. 
    
    Finally, we define $$\Delta_{\Omega^\bullet(P)}: \Omega^\bullet(P) \longrightarrow \Omega^\bullet(P) \otimes H$$ to coincide with $\Delta_P$ on the degree $0$ case, and for degree $1$ elements, $$\Delta_{\Omega^\bullet(P)}(g\,\pi(z))=\Delta(g)(\pi(z)\otimes \mathbbm{1})=g^{(1)}\,\pi(z)\otimes g^{(2)}$$ for all $g$ $\in$ $P$. These constructions provide a differential calculus on $\zeta$. Consequently, $$\Omega^\bullet(B)=B\oplus\Omega^1(B) \qquad \mbox{ with }  \qquad  \Omega^1(B)=\{\lambda\, \pi(z)\mid \lambda \in \C\}.$$ However, this space is not generated by $B$ because $dB=\{ 0\}$.
\end{Example}

\section{The Group Isomorphism}

In this section, we will consider

\begin{Remark}
\label{rema}
Let $\zeta$ be a qpb. Henceforth, we are going to consider that the quantum base space $(B,\cdot,\mathbbm{1},\ast)$ is a  $C^\ast$--algebra (or it can be completed to a $C^\ast$--algebra). Under this assumption and in accordance with \cite{micho3},  for every $\delta^V$ $\in$ $\T$ there exists  $$\{T^\l_k \}^{d_{V}}_{k=1} \subseteq \Mor(\delta^V,\Delta_P)$$ for some $d_{V}$ $\in$ $\N$ such that
\begin{equation}
    \label{generators}
\sum^{d_{V}}_{k=1}x^{V\,\ast}_{ki}x^{V}_{kj}=\delta_{ij}\mathbbm{1},
\end{equation}
with $x^{V}_{ki}:=T^\l_k(e_i)$, where $\T$ is a complete set of mutually non--equivalent irreducible (necessarily finite--dimensional) $\G$--corepresentations with $\delta^\C_\triv$ $\in$ $\T$ (the trivial corepresentation on $\C$) and $\{e_i\}^{n_{V}}_{i=1}$ is the orthonormal basis of $V$ shown in Theorem \ref{rep}.
\end{Remark}

Let $\zeta=(P,B,\Delta_P)$ be a qpb. The surjective map $\beta$ can be used to define the linear isomorphism
\begin{equation*}
\widetilde{\beta}:P\otimes_B P\longrightarrow P\otimes H
\end{equation*}
such that $\widetilde{\beta}(x\otimes_B y)=\beta(x\otimes y)=(x\otimes \mathbbm{1})\cdot \Delta_P(y).$ \cite{stheve}. The degree zero quantum translation map is defined as
\begin{equation}
\qtrs: H \longrightarrow P\otimes_B P,
\end{equation}
such that $\qtrs(g)=\widetilde{\beta}^{-1}(\mathbbm{1}\otimes g).$ Explicitly, by taking the linear basis $\{g^{V}_{ij}\}_{\delta^V,i,j}$  (see Theorem \ref{rep} and Remark \ref{rema}), we have 
\begin{equation}
    \label{qtrs0}
    \qtrs(g^{V}_{ij})=\sum^{d_{V}}_{k=1} x^{V\,\ast}_{ki}\otimes_B x^{V}_{kj}.
\end{equation}
In particular, since $\delta^\C_\triv$  $\in$  $\T$ we have that $\mathbbm{1}$ $\in$ $\{g^{V}_{ij}\}_{\delta^V,i,j}$ and  $$\qtrs(\mathbbm{1})=\mathbbm{1}\otimes_B \mathbbm{1}.$$
We can extend $\qtrs$ to 
\begin{equation}
\label{5.f1.3}
\widetilde{\qtrs}:P\otimes H\longrightarrow P\otimes_B P
\end{equation}
by means of $\widetilde{\qtrs}(x\otimes g^{V}_{ij})=x\,\qtrs(g^{V}_{ij})=\displaystyle\sum^{d_{\alpha}}_{k=1}x\, x^{V\,\ast}_{ki}\otimes_B x^{V}_{kj}.$ A direct calculation shows that  $\widetilde{\beta}$ and $\widetilde{\qtrs}$ are mutually inverse.

Throughout the various computations of this paper, we shall use the symbolic notation
\begin{equation}
\label{5.f1.4}
\qtrs(g)=[g]_1\otimes_B [g]_2.
\end{equation}

Now we shall assume that $\zeta=(P,B,\Delta)$ is endowed with a differential calculus. In this situation $\widetilde{\beta}$ has a natural extension to 
\begin{equation}
\label{6.f1.5}
\widetilde{\beta}:\Omega^\bullet(P)\otimes_{\Omega^\bullet(B)}\Omega^\bullet(P)\longrightarrow \Omega^\bullet(P)\otimes \Gamma^\wedge
\end{equation}
given by $\widetilde{\beta}(w_1\otimes_{\Omega^\bullet(B)}w_2)=(w_1\otimes\mathbbm{1})\cdot \Delta_{\Omega^\bullet(P)}(w_2),$ where the tensor product on the image is the tensor product of graded differential $\ast$--algebras. According to \cite{micho4}, this map is bijective.\\ 

On the other hand, taking a real qpc $\omega$ (which always exists \cite{micho2}) and in accordance with \cite{micho4}, we can extend $\qtrs$ to
\begin{equation}
\label{6.f1.6}
\qtrs:\Gamma\longrightarrow \left(\Omega^\bullet(P)\otimes_{\Omega^\bullet(B)}\Omega^\bullet(P)\right)^{1}
\end{equation}
by means of 
\begin{equation}
\label{6.f1.6.1}
\qtrs(\theta)=\mathbbm{1}\otimes_{\Omega^\bullet(B)}\omega(\theta)-(m_\Omega\otimes_{\Omega^\bullet(B)}\id_{P})(\omega\otimes \qtrs)\ad(\theta)
\end{equation}
when $\theta$ $\in$ $\mathfrak{qg}^\#$, where $$m_\Omega:\Omega^\bullet(P)\otimes\Omega^\bullet(P)\longrightarrow \Omega^\bullet(P)$$ is the product map; and 
\begin{equation}
\label{6.f1.7}
\qtrs(\vartheta\upsilon):= (-1)^{\partial\vartheta\,\partial[\upsilon]_1}[\upsilon]_1\,\qtrs(\vartheta)\,[\upsilon]_2
\end{equation}
    for $\vartheta$ $\in$ $H$, $\upsilon$ $\in$ $\mathfrak{qg}^\#$ or  $\vartheta$ $\in$ $\mathfrak{qg}^\#$, $\upsilon$ $\in$ $H$, with $\qtrs(\upsilon)=[\upsilon]_1\otimes_{\Omega^\bullet(B)} [\upsilon]_2$, where $\partial\vartheta$  is the grade of $\vartheta$. A  direct calculation shows that the corresponding extended map $\widetilde{\qtrs}$ is the right inverse of $\widetilde{\beta}$ for degree $1$ and since $\widetilde{\beta}$ is bijective, we get that $\widetilde{\qtrs}$ is actually its inverse. It is worth mentioning that even when apparently the definition of $\qtrs$ depends on the real qpc $\omega$ chosen, by the uniqueness of the inverse function, the last result tells us that $\qtrs$ is independent of this choice. Since $\widetilde{\beta}$ commutes with the corresponding differential maps, it follows that (\cite{micho4})
 
\begin{Proposition}
\label{6.1.1}
We have  $$\qtrs\circ d=d_{\otimes^\bullet} \circ \qtrs,$$ where $d_{\otimes^\bullet}$ is the differential map of $\Omega^\bullet(P)\otimes_{\Omega^\bullet(B)}\Omega^\bullet(P)$.
\end{Proposition}

Now let us take the universal graded differential $\ast$ --algebra $(\otimes^\bullet_H\Gamma,d_{\otimes_H},\ast)$ of  $(\Gamma,d)$ \cite{stheve}. The quantum translation map can be extended naturally to $\otimes^\bullet_H\Gamma$ by means of
\begin{equation}
\label{6.f1.8}
\qtrs(\vartheta\otimes_H \upsilon):= (-1)^{\partial\vartheta\,\partial[\upsilon]_1} \,[\upsilon]_1\,\qtrs(\vartheta)\,[\upsilon]_2
\end{equation} 
if $\qtrs(\upsilon)= [\upsilon]_1\otimes_{\Omega^\bullet(B)}[\upsilon]_2$. By induction and Proposition \ref{6.1.1} it can be proved that
\begin{equation}
\label{6.f1.9}
\qtrs(dg_0\otimes_H dg_1\otimes_H...\otimes_H dg_k)=d_{\otimes^\bullet}(\qtrs(g_0\,dg_1\otimes_H...\otimes_H dg_k))
\end{equation}
for $g_1$,..., $g_k$ $\in$ $H$.  Now by considering the bilateral ideal $\mathcal{Q}$ generated by the relation presented on equation (\ref{udtensor1}), a direct calculation shows that $\qtrs(\mathcal{Q})=0$ and therefore,  $\qtrs$ can be defined on the universal differential envelope $\ast$--calculus 
\begin{equation}
\label{6.f1.10}
\qtrs:\Gamma^\wedge\longrightarrow \Omega^\bullet(P)\otimes_{\Omega^\bullet(B)}\Omega^\bullet(P).
\end{equation}
As before, the corresponding extended map
\begin{equation}
\label{6.f1.12}
\widetilde{\qtrs}:\Omega^\bullet(P)\otimes \Gamma^\wedge\longrightarrow \Omega^\bullet(P)\otimes_{\Omega^\bullet(B)}\Omega^\bullet(P).
\end{equation}
given by $\widetilde{\qtrs}(w_1\otimes \vartheta)=(w_1\otimes \mathbbm{1})\cdot \qtrs(\vartheta) $ is the inverse map of $\widetilde{\beta}$ \cite{micho4}. According to \cite{micho5} we have

\begin{Proposition}
\label{qtrsprop}
The following properties hold
   \begin{enumerate}
\item $[\vartheta]_1\,[\vartheta]_2=\epsilon(\vartheta)\mathbbm{1}$.
\item $(\id_{\Omega^\bullet(P)}\otimes_{\Omega^\bullet(B)}\Delta_{\Omega^\bullet(P)})\circ\qtrs=(\qtrs\otimes \id_{\Gamma^\wedge})\circ \Delta$.
\item $(\Delta_{\Omega^\bullet(P)}\otimes_{\Omega^\bullet(B)} \id_{\Omega^\bullet(P)})\circ\qtrs=(\sigma\otimes_{\Omega^\bullet(B)} \id_{\Omega^\bullet(P)}) \circ(S\otimes \qtrs)\circ \Delta$, where $$\sigma:\Gamma^\wedge \otimes \Omega^\bullet(P) \longrightarrow \Omega^\bullet(P)\otimes\Gamma^\wedge $$ is the canonical graded twist map, i.e., $\sigma(\vartheta\otimes w)=(-1)^{kl}\,w\otimes \vartheta$ if $w$ $\in$ $\Omega^k(P)$ and $\vartheta$ $\in$ $\Gamma^{\wedge l}$.
\item $\mu\,\qtrs(\vartheta)=(-1)^{lk}\qtrs(\vartheta)\,\mu$ for all $\mu$ $\in$ $\Omega^k(B)$, $\vartheta$ $\in$ $\Gamma^{\wedge l}.$
\end{enumerate} 
\end{Proposition}

Let $$\f_1,\,\f_2:\Gamma^\wedge\longrightarrow \Omega^\bullet(P) $$ be two graded linear maps. The convolution product of $\f_1$ with $\f_2$ is defined by $$\f_1\ast \f_2=m_\Omega \circ (\f_1\otimes \f_2)\circ \Delta:\Gamma^\wedge \longrightarrow \Omega^\bullet(P).$$

Henceforth we will just consider graded maps $\f$ such that 
\begin{equation}
    \label{4.f3}
    \f(\mathbbm{1})=\mathbbm{1} \qquad \mbox{and} \qquad (\f\otimes \id_{\Gamma^\wedge})\circ \Ad=\Delta_{\Omega^\bullet(P)}\circ \f,
\end{equation}
where $\Ad: \Gamma^\wedge\longrightarrow \Gamma^\wedge \otimes \Gamma^\wedge$ is the extension of the right $\G$--coaction $\Ad: H\longrightarrow H\otimes H$ (see equation (\ref{2.f11})). We say that $\f$ is a {\it convolution invertible map} if there exists a graded linear map $$\f^{-1}: \Gamma^\wedge\longrightarrow \Omega^\bullet(P)$$
such that 
\begin{equation}
    \label{4.f3.1}
    \f\ast \f^{-1}=\f^{-1}\ast \f =\mathbbm{1}\epsilon.
\end{equation}
A direct calculation shows that the set of all convolution invertible maps $\{ \f: \Gamma^\wedge\longrightarrow \Omega^\bullet(P)\} $ is a group with respect to the convolution product. The next proposition is one of the purposes of this paper

\begin{Proposition}
    \label{4.1}
    There exist a group isomorphism between the group of all convolution invertible maps and the group of all graded left $\Omega^\bullet(B)$--module isomorphisms $$\F:\Omega^\bullet(P)\longrightarrow \Omega^\bullet(P) $$ that satisfy 
    \begin{equation}
        \label{4.f4}
        \F(\mathbbm{1})=\mathbbm{1} \qquad \mbox{ and }\qquad (\F\otimes \id_{\Gamma^\wedge})\circ \Delta_{\Omega^\bullet(P)}=\Delta_{\Omega^\bullet(P)}\circ \F.
    \end{equation}
    Here we are considering the group product $(\F_1\circ \F_2)(w)=\F_2(\F_1(w))$.
\end{Proposition}

\begin{proof}
Let us start by considering a map $\F$. Then we define a graded linear map 
\begin{equation}
\label{4.f5}
\f_\F:=m_{\Omega^\bullet}\circ (\id_{\Omega^\bullet(P)}\otimes_{\Omega^\bullet(B)}\F)\circ \qtrs:\Gamma^\wedge\longrightarrow \Omega^\bullet(P),
\end{equation}
where $m_{\Omega^\bullet}:\Omega^\bullet(P)\otimes_{\Omega^\bullet(B)}\Omega^\bullet(P)\longrightarrow \Omega^\bullet(P)$ is the product map. We are going to show that $\f_\F$ is a convolution invertible map. First, since $\qtrs(\mathbbm{1})=\mathbbm{1}\otimes_B \mathbbm{1}$ it follows $\f_\F(\mathbbm{1})=\mathbbm{1}.$ Secondly, according to Proposition \ref{qtrsprop} point $3$
\begin{eqnarray}
\label{9.f1.16}
\Delta_{\Omega^\bullet (P)}\circ \f_\F &= & \widehat{m}_{\Omega^\bullet}\circ (\Delta_{\Omega^\bullet (P)}\otimes_{\Omega^\bullet(B)}(\Delta_{\Omega^\bullet (P)} \circ \F))\circ \qtrs \nonumber
  \\
  &= &
\widehat{m}_{\Omega^\bullet}\circ (\id_{\Omega^\bullet(P)\otimes \Gamma^\wedge}\otimes_{\Omega^\bullet(B)}(\Delta_{\Omega^\bullet (P)} \circ \F))\circ (\Delta_{\Omega^\bullet (P)}\otimes_{\Omega^\bullet(B)}\id_{\Omega^\bullet(P)})\circ \qtrs 
  \\
  &= &
\widehat{m}_{\Omega^\bullet}\circ (\id_{\Omega^\bullet(P)\otimes \Gamma^\wedge}\otimes_{\Omega^\bullet(B)}(\Delta_{\Omega^\bullet (P)} \circ \F))\circ (\sigma \otimes_{\Omega^\bullet(B)}\id_{\Omega^\bullet(P)})\circ (S\otimes \qtrs) \circ \Delta,
\nonumber
\end{eqnarray}
where $$\widehat{m}_{\Omega^\bullet}:(\Omega^\bullet(P)\otimes \Gamma^\wedge)\otimes_{\Omega^\bullet(B)} (\Omega^\bullet(P)\otimes \Gamma^\wedge) \longrightarrow \Omega^\bullet(P)\otimes \Gamma^\wedge$$ is such that $\widehat{m}_{\Omega^\bullet}(w_1\otimes \vartheta_1\otimes_{\Omega^\bullet(B)}w_2\otimes \vartheta_2)=(-1)^{kl}w_1w_2\otimes \vartheta_1\vartheta_2$ if $w_2$ $\in$ $\Omega^k(P)$ and $\vartheta_1$ $\in$ $\Gamma^{\wedge l}$. On the other hand, by equation (\ref{4.f4}) and Proposition \ref{qtrsprop} point $2$ 
\begin{eqnarray}
\label{9.f1.17}
\Delta_{\Omega^\bullet (P)}\circ \f_\F &= & \widehat{m}_{\Omega^\bullet}\circ (\Delta_{\Omega^\bullet (P)}\otimes_{\Omega^\bullet(B)}(\Delta_{\Omega^\bullet (P)} \circ \F))\circ \qtrs \nonumber
  \\
  &= &
\widehat{m}_{\Omega^\bullet}\circ (\Delta_{\Omega^\bullet (P)}\otimes_{\Omega^\bullet(B)}((\F\otimes\id_{\Gamma^\wedge}) \circ \Delta_{\Omega^\bullet (P)}))\circ \qtrs
  \\
  &= &
\widehat{m}_{\Omega^\bullet}\circ (\Delta_{\Omega^\bullet (P)}\otimes_{\Omega^\bullet(B)}(\F\otimes\id_{\Gamma^\wedge}))\circ (\id_{\Omega^\bullet(P)}\otimes_{\Omega^\bullet(B)}\Delta_{\Omega^\bullet (P)}) \circ \qtrs \nonumber
  \\
  &= &
\widehat{m}_{\Omega^\bullet}\circ (\Delta_{\Omega^\bullet (P)}\otimes_{\Omega^\bullet(B)}(\F\otimes\id_{\Gamma^\wedge}))\circ (\qtrs\otimes \id_{\Gamma^\wedge}) \circ \Delta; \nonumber
\end{eqnarray}
but by considering again Proposition \ref{qtrsprop} point $3$
\begin{eqnarray*}
\Delta_{\Omega^\bullet (P)}\circ \f_\F &= & \widehat{m}_{\Omega^\bullet}\circ (\Delta_{\Omega^\bullet (P)}\otimes_{\Omega^\bullet(B)}(\F\otimes\id_{\Gamma^\wedge}))\circ (\qtrs\otimes \id_{\Gamma^\wedge}) \circ \Delta 
  \\
  &= &
\widehat{m}_{\Omega^\bullet}\circ (\id_{\Omega^\bullet(P)\otimes \Gamma^\wedge}\otimes_{\Omega^\bullet(B)}(\F\otimes\id_{\Gamma^\wedge}))\circ (\Delta_{\Omega^\bullet (P)}\otimes_{\Omega^\bullet(B)}\id_{\Omega^\bullet(P)\otimes \Gamma^\wedge})\nonumber
  \\
  & &
\circ\,(\qtrs\otimes \id_{\Gamma^\wedge}) \circ \Delta \nonumber
  \\
  &= &
\widehat{m}_{\Omega^\bullet}\circ (\id_{\Omega^\bullet(P)\otimes \Gamma^\wedge}\otimes_{\Omega^\bullet(B)}(\F\otimes\id_{\Gamma^\wedge}))\nonumber
  \\
  & &
\circ\,[(({\sigma}\otimes_{\Omega^\bullet(B)}\id_{\Omega^\bullet(P)})\circ (S\otimes \qtrs)\circ \Delta)\otimes \id_{\Gamma^\wedge}]\circ \Delta. \nonumber
\end{eqnarray*}
So for all $\vartheta$ $\in$ $\Gamma^\wedge$ 
\begin{eqnarray*}
\Delta_{\Omega^\bullet (P)}(\f_\F(\vartheta)) &= & (-1)^{\partial\vartheta^{(1)}(\partial[\vartheta^{(2)}]_1+\partial[\vartheta^{(2)}]_2)}\,[\vartheta^{(2)}]_1F([\vartheta^{(2)}]_2)\otimes S(\vartheta^{(1)})\vartheta^{(3)} 
  \\
  &= &
(-1)^{\partial\vartheta^{(1)}\partial\vartheta^{(2)}}\,[\vartheta^{(2)}]_1F([\vartheta^{(2)}]_2)\otimes S(\vartheta^{(1)})\vartheta^{(3)}\;=\;(\f_\F\otimes\id_{\Gamma^\wedge})\Ad(\vartheta),
\end{eqnarray*}
\noindent because $\qtrs$ is a graded linear map. Finally, consider $$\f_{\F^{-1}}:=m_{\Omega^\bullet}\circ (\id_{\Omega^\bullet(P)}\otimes_{\Omega^\bullet(B)}\F^{-1})\circ \qtrs:\Gamma^\wedge\longrightarrow \Omega^\bullet(P).$$   Then for all $\vartheta$ $\in$ $\Gamma^\wedge$ 
\begin{equation*}
    (\f_\F\ast \f_{\F^{-1}})(\vartheta)=[\vartheta^{(1)}]_1\,\underbrace{\F([\vartheta^{(1)}]_2)\,[\vartheta^{(2)}]_1}\,\F^{-1}([\vartheta^{(2)}]_2).
\end{equation*}
We claim that the expression in the brace is an element of $\Omega^\bullet(B)$. Indeed, notice that
\begin{eqnarray}
    \label{ec.99}
\Delta_{\Omega^\bullet (P)}((\f_\F\ast \f_{\F^{-1}})(\vartheta))&=&\Delta_{\Omega^\bullet (P)}\left([\vartheta^{(1)}]_1\,\F([\vartheta^{(1)}]_2)\,[\vartheta^{(2)}]_1\,\F^{-1}([\vartheta^{(2)}]_2)\right)\nonumber
\\
 &=&
 \Delta_{\Omega^\bullet (P)}([\vartheta^{(1)}]_1)\,\Delta_{\Omega^\bullet (P)}(\F([\vartheta^{(1)}]_2)\,[\vartheta^{(2)}]_1)\,\Delta_{\Omega^\bullet (P)}(\F^{-1}([\vartheta^{(2)}]_2)).
\end{eqnarray}
Also we know 
$$\Delta_{\Omega^\bullet (P)}((\f_\F\ast \f_{\F^{-1}})(\vartheta))=\Delta_{\Omega^\bullet (P)}(\f_\F(\vartheta^{(1)}))\cdot\Delta_{\Omega^\bullet (P)}(\f_{\F^{-1}}(\vartheta^{(2)})).$$ By applying equation (\ref{9.f1.17}) on $\Delta_{\Omega^\bullet (P)}(\f_\F(\vartheta^{(1)}))$  and equation (\ref{9.f1.16}) on $\Delta_{\Omega^\bullet (P)}(\f_\F(\vartheta^{(2)}))$ we have
\begin{eqnarray*}
\Delta_{\Omega^\bullet (P)}((\f_\F\ast \f_{\F^{-1}})(\vartheta))  &= &
(-1)^{\partial\vartheta^{(3)}\partial[\vartheta^{(4)}]_1}\,\Delta_{\Omega^\bullet (P)}([\vartheta^{(1)}]_1)\,(\F([\vartheta^{(1)}]_2)\otimes \vartheta^{(2)})
  \\
 & &
([\vartheta^{(4)}]_1\otimes S(\vartheta^{(3)}))\,\Delta_{\Omega^\bullet (P)}(\F^{-1}([\vartheta^{(4)}]_2))
  \\
  &= &
(-1)^{(\partial\vartheta^{(3)}+\partial\vartheta^{(2)})\partial[\vartheta^{(4)}]_1}\, \Delta_{\Omega^\bullet (P)}([\vartheta^{(1)}]_1)\, (\F([\vartheta^{(1)}]_2)[\vartheta^{(4)}]_1\otimes \vartheta^{(2)}S(\vartheta^{(3)}))
  \\
  & &
\Delta_{\Omega^\bullet (P)}(\F^{-1}([\vartheta^{(4)}]_2)).
\end{eqnarray*}
On the other hand, by the definition of $\epsilon$ (see equation (\ref{2.f3.6})) and the graded Hopf $\ast$--algebra structure of $\Gamma^{\wedge\,\infty}$ and we get
$$\Delta_{\Omega^\bullet (P)}([\vartheta^{(1)}]_1)\, (\F([\vartheta^{(1)}]_2)[\vartheta^{(2)}]_1\otimes \mathbbm{1})\,\Delta_{\Omega^\bullet (P)}(\F^{-1}([\vartheta^{(2)}]_2)) $$  $$=\Delta_{\Omega^\bullet (P)}([\vartheta^{(1)}]_1)\, (\F([\vartheta^{(1)}]_2)\otimes \mathbbm{1})
 \widetilde{\beta}((\id_{\Omega^\bullet(P)}\otimes_{\Omega^\bullet(B)}\F^{-1})\qtrs(\vartheta^{(2)}))
  $$ $$
 =\Delta_{\Omega^\bullet (P)}([\vartheta^{(1)}]_1)\, (\F([\vartheta^{(1)}]_2)\otimes \hspace{-0.7cm}\underbrace{\epsilon(\vartheta^{(2)})}_{\not=0 \mathrm{\, only\, for\, } \Gamma^{\wedge\,0}=H}\hspace{-0.7cm}\mathbbm{1})\,
 \widetilde{\beta}((\id_{\Omega^\bullet(P)}\otimes_{\Omega^\bullet(B)}\F^{-1})\qtrs(\vartheta^{(3)}))
 $$  $$= (-1)^{\partial\vartheta^{(2)}\partial[\vartheta^{(4)}]_1}  \, \Delta_{\Omega^\bullet (P)}([\vartheta^{(1)}]_1)\, (\;\F([\vartheta^{(1)}]_2)\otimes \underbrace{\epsilon(\vartheta^{(2)})\mathbbm{1}}_{0-\mathrm{degree}})\,\widetilde{\beta}((\id_{\Omega^\bullet(P)}\otimes_{\Omega^\bullet(B)}\F^{-1})\qtrs(\vartheta^{(3)})).$$ $$ =(-1)^{(\partial\vartheta^{(3)}+\partial\vartheta^{(2)})\partial[\vartheta^{(4)}]_1}\, \Delta_{\Omega^\bullet (P)}([\vartheta^{(1)}]_1)\, (\F([\vartheta^{(1)}]_2)\otimes \underbrace{\vartheta^{(2)}S(\vartheta^{(3)}}_{0-\mathrm{degree}}))\,\widetilde{\beta}((\id_{\Omega^\bullet(P)}\otimes_{\Omega^\bullet(B)}\F^{-1})\qtrs(\vartheta^{(4)}))$$
$$=(-1)^{(\partial\vartheta^{(3)}+\partial\vartheta^{(2)})\partial[\vartheta^{(4)}]_1}\, \Delta_{\Omega^\bullet (P)}([\vartheta^{(1)}]_1)\, (\F([\vartheta^{(1)}]_2)\otimes \underbrace{\vartheta^{(2)}S(\vartheta^{(3)}}_{0-\mathrm{degree}}))\, ([\vartheta^{(4)}]_1\otimes \mathbbm{1})\,\Delta_{\Omega^\bullet (P)}(\F^{-1}([\vartheta^{(4)}]_2))$$
$$=(-1)^{(\partial\vartheta^{(3)}+\partial\vartheta^{(2)})\partial[\vartheta^{(4)}]_1}\, \Delta_{\Omega^\bullet (P)}([\vartheta^{(1)}]_1)\, (\F([\vartheta^{(1)}]_2)[\vartheta^{(4)}]_1\otimes \vartheta^{(2)}S(\vartheta^{(3)}))\Delta_{\Omega^\bullet (P)}(\F^{-1}([\vartheta^{(4)}]_2))$$ and hence $$\Delta_{\Omega^\bullet (P)}((\f_\F\ast \f_{\F^{-1}})(\vartheta))=\Delta_{\Omega^\bullet (P)}([\vartheta^{(1)}]_1)\, (\F([\vartheta^{(1)}]_2)[\vartheta^{(2)}]_1\otimes \mathbbm{1})\,\Delta_{\Omega^\bullet (P)}(\F^{-1}([\vartheta^{(2)}]_2)).$$ By equation (\ref{ec.99}) and the last equality we conclude $$\Delta_{\Omega^\bullet (P)}([\vartheta^{(1)}]_1)\,\Delta_{\Omega^\bullet (P)}(\F([\vartheta^{(1)}]_2)\,[\vartheta^{(2)}]_1)\,\Delta_{\Omega^\bullet (P)}(\F^{-1}([\vartheta^{(2)}]_2))$$ $$=\Delta_{\Omega^\bullet (P)}([\vartheta^{(1)}]_1)\, (\F([\vartheta^{(1)}]_2)[\vartheta^{(2)}]_1\otimes \mathbbm{1})\,\Delta_{\Omega^\bullet (P)}(\F^{-1}([\vartheta^{(2)}]_2)),$$  which  proves our claim. Since $\F^{-1}$ is a left $\Omega^\bullet(B)$--module morphism, by our previous claim and Proposition \ref{qtrsprop} point $1$ 
\begin{eqnarray*}
(\f_\F\ast \f_{\F^{-1}})(\vartheta)\;=\; [\vartheta^{(1)}]_1\,\F([\vartheta^{(1)}]_2)\,[\vartheta^{(2)}]_1\,\F^{-1
}([\vartheta^{(2)}]_2) &=& [\vartheta^{(1)}]_1\,\F^{-1}(\F([\vartheta^{(1)}]_2)\,[\vartheta^{(2)}]_1[\vartheta^{(2)}]_2)
  \\
 &=&
 [\vartheta^{(1)}]_1\,\F^{-1}(\F([\vartheta^{(1)}]_2)\,\epsilon(\vartheta^{(2)}))
   \\
 &=&
\epsilon(\vartheta^{(1)})\epsilon(\vartheta^{(2)})\mathbbm{1} 
   \\
 &=&
 \epsilon(\vartheta)\mathbbm{1}
\end{eqnarray*}
for every $\vartheta$ $\in$ $\Gamma^\wedge$. In a similar way, we can prove that $\f_{\F^{-1}}\ast \f_\F=\mathbbm{1}\epsilon$ and hence $\f_\F$ is a convolution invertible map.\\

Conversely, for a given convolution invertible map $\f$ let us define the graded linear map
\begin{equation}
\label{4.f5.56}
\F_\f:=m_{\Omega}\circ (\id_{\Omega^\bullet(P)}\otimes \f)\circ \Delta_{\Omega^\bullet (P)}:\Omega^\bullet(P)\longrightarrow \Omega^\bullet(P).
\end{equation}
We are going to prove that $\F_\f$ is a graded left $\Omega^\bullet(B)$--module isomorphism which satisfies equation (\ref{4.f4}). First of all, by equation (\ref{4.f3}) it is obvious that $\F_\f(\mathbbm{1})=\mathbbm{1}$. Secondly, taking $\mu$ $\in$ $\Omega^\bullet(B)$ and $w$ $\in$ $\Omega^\bullet(P)$ we have $\F_\f(\mu w)=m_{\Omega}(\id_{\Omega^\bullet(P)}\otimes \f) (\mu\otimes \mathbbm{1}) \Delta_{\Omega^\bullet (P)}(w)=\mu\, \F_f(w).$ Thirdly, by equation (\ref{4.f3}) 
\begin{eqnarray*}
 \Delta_{\Omega^\bullet (P)}\circ \F_\f &=& \widehat{m}_{\Omega}\circ (\Delta_{\Omega^\bullet (P)}\otimes(\Delta_{\Omega^\bullet (P)}\circ \f))\circ \Delta_{\Omega^\bullet (P)} 
  \\
 &=&
\widehat{m}_{\Omega}\circ(\Delta_{\Omega^\bullet (P)} \otimes((\f\otimes\id_{\Gamma^\wedge})\circ \Ad))\circ \Delta_{\Omega^\bullet (P)}, 
\end{eqnarray*}
where $$\widehat{m}_{\Omega}:(\Omega^\bullet(P)\otimes \Gamma^\wedge)\otimes (\Omega^\bullet(P)\otimes \Gamma^\wedge) \longrightarrow \Omega^\bullet(P)\otimes \Gamma^\wedge$$ is such that $\widehat{m}_{\Omega^\bullet}(w_1\otimes \vartheta_1\otimes w_2\otimes \vartheta_2)=(-1)^{kl}w_1w_2\otimes \vartheta_1\vartheta_2$ if $w_2$ $\in$ $\Omega^k(P)$ and $\vartheta_1$ $\in$ $\Gamma^{\wedge l}$. In this way, by the graded Hopf $\ast$--algebra structure of $\Gamma^{\wedge\,\infty}$ and the fact that $\epsilon\not=0$ only for $\Gamma^{\wedge\,0}=H$, we get for all $w$ $\in$ $\Omega^\bullet(P)$ 
\begin{eqnarray*}
 \Delta_{\Omega^\bullet (P)}(\F_\f(w)) &=& (-1)^{\partial w^3(\partial w^1+\partial w^2)} w^{(0)}\f(w^{(3)})\otimes w^{(1)}S(w^{(2)})w^{(4)} 
  \\
 &=&
(-1)^{\partial w^3\partial w^1} w^{(0)}\f(w^{(2)})\otimes \epsilon(w^{(1)})w^{(3)}
\\
 &=&
 w^{(0)}\f(w^{(2)})\otimes \epsilon(w^{(1)})w^{(3)}
 \\
 &=&
 w^{(0)}\f(w^{(2)})\otimes \epsilon(w^{(1)})w^{(3)}
 \\
 &=&
 w^{(0)}\f(w^{(1)})\otimes w^{(2)}=(\F_\f\otimes \id_{\Gamma^\wedge})\Delta_{\Omega^\bullet (P)}(w).
\end{eqnarray*}
Finally, notice that for all $w$ $\in$ $\Omega^\bullet(P)$
\begin{eqnarray*}
\F_{\f^{-1}}(\F_{\f}(w))  \;=\;  m_{\Omega}(\id_{\Omega^\bullet(P)}\otimes \f^{-1})\,\Delta_{\Omega^\bullet (P)}(\F_\f(w)) &=& m_{\Omega}(\id_{\Omega^\bullet(P)}\otimes \f^{-1})(\F_\f\otimes \id_{\Gamma^\wedge})\,\Delta_{\Omega^\bullet (P)}(w)
 \\
 &=&
\F_f(w^{(0)})\,\f^{-1}(w^{(1)})
 \\
 &=&
w^{(0)}\,\f(w^{(1)})\,\f^{-1}(w^{(2)})
 \\
 &=&
w^{(0)}\,\epsilon(w^{(1)})\;=\;w.
\end{eqnarray*}
A similar calculation shows $\F_\f(\F_{\f^{-1}}(w))=w$, so $\F_\f$ is biyective and $\F^{-1}_\f=\F_{\f^{-1}}.$\\ 

Our next step is to prove that $$\F \xmapsto{\widehat{\Lambda}} \f_\F\,,\qquad \f\xmapsto{\widetilde{\Lambda}} \F_\f $$ are mutually inverse. Notice that for all $w$ $\in$ $\Omega^\bullet(P)$ $$\F_{\f_\F}(w)=m_{\Omega}(\id_{\Omega^\bullet(P)}\otimes (m_{\Omega^\bullet}(\id_{\Omega^\bullet(P)}\otimes_{\Omega^\bullet(B)}\F)\qtrs))\Delta_{\Omega^\bullet (P)}(w)=\underbrace{w^{(0)}[w^{(1)}]_1}\F([w^{(1)}]_2),$$
where the expression in the brace is an element of $\Omega^\bullet(B)$. In fact, 
$$\Delta_{\Omega^\bullet(P)}(\F_{\f_\F}(w))=\Delta_{\Omega^\bullet(P)}\left(w^{(0)}[w^{(1)}]_1\F([w^{(1)}]_2)\right)=\Delta_{\Omega^\bullet(P)}(w^{(0)}[w^{(1)}]_1)\, \Delta_{\Omega^\bullet(P)}(\F([w^{(1)}]_2));$$
however by equation (\ref{9.f1.16}) we get

\begin{eqnarray*}
\Delta_{\Omega^\bullet (P)}(\F_{\f_\F}(w)) \;=\; \Delta_{\Omega^\bullet (P)}(w^{(0)}\f_\F(w^{(1)})) &= & (w^{(0)}\otimes w^{(1)})\,\Delta_{\Omega^\bullet (P)}(\f_\F(w^{(2)}))
  \\
  &= &
(-1)^{\partial w^{(2)}\partial[w^{(3)}]_1}(w^{(0)}\otimes w^{(1)})
  \\
  & &
([w^{(3)}]_1\otimes S(w^{(2)}))\,\Delta_{\Omega^\bullet (P)}(\F([w^{(3)}]_2))
  \\
  &= &
(-1)^{(\partial w^{(1)}+\partial w^{(2)})\partial[w^{(3)}]_1}\, (w^{(0)}[w^{(3)}]_1\otimes w^{(1)}S(w^{(2)}))
  \\
  & &
\Delta_{\Omega^\bullet (P)}(\F([w^{(3)}]_2))
\\
  &= &
  (-1)^{\partial w^{(1)}\partial[w^{(2)}]_1}\, (w^{(0)}[w^{(2)}]_1\otimes \hspace{-0.6cm}\underbrace{\epsilon(w^{(1)})}_{\not=0 \mathrm{\, only\, for \, \Gamma^{\wedge\,0}=H}}\hspace{-0.6cm} \mathbbm{1})
  \\
  & & \Delta_{\Omega^\bullet (P)}(\F([w^{(2)}]_2))
   \\
  &= &
  (w^{(0)}[w^{(2)}]_1\otimes \epsilon(w^{(1)})\mathbbm{1})\,\Delta_{\Omega^\bullet (P)}(\F([w^{(2)}]_2))
\\
  &= &
(w^{(0)}\otimes \epsilon(w^{(1)})\mathbbm{1})([w^{(2)}]_1\otimes \mathbbm{1})\,\Delta_{\Omega^\bullet (P)}(\F([w^{(2)}]_2))
    \\
 &=&
 (w^{(0)}\otimes \epsilon(w^{(1)})\mathbbm{1})\,\widetilde{\beta}((\id_{\Omega^\bullet(P)}\otimes_{\Omega^\bullet(B)}\F)\qtrs(w^{(2)}))
  \\
 &=&
(w^{(0)}\otimes \mathbbm{1})\,\widetilde{\beta}((\id_{\Omega^\bullet(P)}\otimes_{\Omega^\bullet(B)}\F)\qtrs(w^{(1)}))
  \\
 &=& 
(w^{(0)}[w^{(1)}]_1\otimes \mathbbm{1})\, \Delta_{\Omega^\bullet (P)}(\F([w^{(1)}]_2)). 
\end{eqnarray*}
This implies that $$\Delta_{\Omega^\bullet(P)}(w^{(0)}[w^{(1)}]_1)\,  \Delta_{\Omega^\bullet(P)}(\F([w^{(1)}]_2))= (w^{(0)}[w^{(1)}]_1\otimes \mathbbm{1})\, \Delta_{\Omega^\bullet (P)}(\F([w^{(1)}]_2)),$$ which proves our  assertion. Since $\F$ is a left $\Omega^\bullet(B)$--module morphism, by this assertion and Proposition \ref{qtrsprop} point $1$, we have  $$\F_{\f_\F}(w)=w^{(0)}[w^{(1)}]_1\,\F([w^{(1)}]_2)=\F(w^{(0)}[w^{(1)}]_1[w^{(1)}]_2))=\F(w^{(0)}\epsilon(w^{(1)}))=\F(w).$$  

On the other hand, for every $\vartheta$ $\in$ $\Gamma^\wedge$

\begin{eqnarray*}
\f_{\F_\f}(\vartheta)  &=&  m_{\Omega^\bullet}(\id_{\Omega^\bullet(P)}\otimes_{\Omega^\bullet(B)}(m_{\Omega}(\id_{\Omega^\bullet(P)}\otimes \f)\Delta_{\Omega^\bullet (P)}))\qtrs(\vartheta)
 \\
 &=&
[\vartheta]_1[\vartheta]^{(0)}_2\f([\vartheta]^{(1)}_2)
 \\
 &=&
m_{\Omega}(\id_{\Omega^\bullet(P)}\otimes \f)([\vartheta]_1\otimes\mathbbm{1})([\vartheta]^{(0)}_2\otimes [\vartheta]^{(1)}_2)
\\
 &=&
m_{\Omega}(\id_{\Omega^\bullet(P)}\otimes \f)([\vartheta]_1\otimes\mathbbm{1})\Delta_{\Omega^\bullet(P)}([\vartheta]_2)
 \\
 &=&
m_{\Omega}(\id_{\Omega^\bullet(P)}\otimes \f)\,\widetilde{\beta}(\qtrs(\vartheta))
 \\
 &=&
m_{\Omega}(\id_{\Omega^\bullet(P)}\otimes \f)(\mathbbm{1}\otimes \vartheta)\;=\;\f(\vartheta)
\end{eqnarray*}

\noindent and hence $\widetilde{\Lambda}=\widehat{\Lambda}^{-1}$. Finally, to complete the proof it is enough to prove that $\widehat{\Lambda}$ or $\widehat{\Lambda}^{-1}$ is a group morphism. In fact $$\F_{\f_1\ast \f_2}(w) =  m_\Omega(\id_{\Omega^\bullet(P)}\otimes (\f_1\ast \f_2))\Delta_{\Omega^\bullet (P)}(w)= w^{(0)}\,\f_1(w^{(1)})\,\f_2(w^{(2)});$$ while by equation (\ref{4.f4})
\begin{eqnarray*}
\F_{\f_2}(\F_{\f_1}(w))  \;=\;  m_{\Omega}(\id_{\Omega^\bullet(P)}\otimes \f_2)\,\Delta_{\Omega^\bullet (P)}(\F_{\f_1}(w)) &=& m_{\Omega}(\id_{\Omega^\bullet(P)}\otimes \f_2)(\F_{\f_1}\otimes \id_{\Gamma^\wedge})\,\Delta_{\Omega^\bullet (P)}(w)
 \\
 &=&
\F_{\f_1}(w^{(0)})\,\f_2(w^{(1)})
 \\
 &=&
w^{(0)}\,\f_1(w^{(1)})\,\f_2(w^{(2)})
\end{eqnarray*}
for all $w$ $\in$ $\Omega^\bullet(P)$. Therefore $\f_1\ast \f_2\xmapsto{\widehat{\Lambda}^{-1}} \F_{\f_1} \circ \F_{\f_2}$ and the proposition follows.
\end{proof}

\section{The form of the maps $\{T^\l_k \}$ and the Gauge Principle in Differential Geometry}

Let $G$ be a compact matrix Lie group, and let $\G$ be its associated quantum group \cite{woro1}. If $\delta^V$ is a (unitary) irreducible $\G$--corepresentation, then it induces a $G$--representation on $V$. Indeed, by Theorem \ref{rep} we know that $$\delta^V(e_j)=\sum^n_{i=1}e_i\otimes g^V_{ij}.$$ Then the linear map $$\alpha^V: G\times V\longrightarrow V $$ given by $$\alpha^V(A,e_j)=\sum^n_{i=1}g^V_{ij}(A)\,e_i$$ is a unitary and irreducible $G$--representation. The following proposition is another of our purposes.

\begin{Proposition}
    Let $G$ be a compact matrix Lie group and let $\pi:P\longrightarrow B$ be a classical principal $G$--bundle, where $P$ is the total space, $B$ is the base space and $\pi$ is the bundle projection. Assume $B$ are compact. If $\T$ is a complete set of mutually non--equivalent irreducible $\G$--corepresentations with $\delta^\C_\triv$ $\in$ $\T$, then for every  $\delta^V$ $\in$ $\T$, there exists $\{T^\l_k\}^{d_V}_{k=1}$ $\subseteq$ $\Mor(\delta^V, \Delta_P)$, for some $d_V$ $\in$ $\N$, such that equation (\ref{generators}) holds. Here $$\Delta_P:C^\infty_\C(P)\longrightarrow C^\infty_\C(P\times G)\cong C^\infty_\C(P)\otimes C^\infty_\C(G)$$ is the pull--back of the right $G$--action on $P$. In this setting, $C^\infty_\C(P)$ denotes the space of $\C$--valued smooth functions of $P$ and the tensor product is taken to be the completed injective tensor product.
\end{Proposition}

\begin{proof}
     Let $\delta^V$ $\in$ $\T$ and consider its associated representation $\alpha^V$ with $n=\mathrm{dim}_\C(V)$. Since $\alpha^V$ is unitary, the associated vector bundle $\pi_{\alpha^V}: E^V \longrightarrow B$ is a Hermitian bundle. 
     
     For each $b$ $\in$ $B$, let $(U_b, \Phi_b)$ be a local trivialization of the Hermitian bundle $\pi_{\alpha^V}: E^V \longrightarrow B$ around $b$ associated to  a principal $G$--bundle local trivialization $(U_b, \Psi_b)$ of $\pi:P\longrightarrow B$. Then there exists a set of local sections $\{\hat{s}^{b}_1, \ldots, \hat{s}^{b}_1\}$ $\subseteq$ $\Gamma\bigl(\pi_{{\alpha^V}}^{-1}(U_b)\bigr)$ such that $\hat{s}^b_i(a)=(a,e_i)$ for all $a$ $\in$ $U_b$ (under the diffeomorphism $\Phi_b$). Since $\{U_b \}_{b\in B}$ is an open cover, by compactness, there exist points $b_1, \ldots, b_r \in B$ such that $\{\,U_{b_i}\}_{i=1}^r $ remains an open cover of $B$. Let $\{\rho_{b_i} \}_{i=1}^r$ be a partition of unity subordinate to the open cover $\{\,U_{b_i}\}_{i=1}^r$, where each $\rho_{b_i}$ has compact support and admits a smooth square root.  In this way, consider the global sections of $E^V$ $$\{s^i_j=\sqrt{\rho_{b_i}} \,\hat{s}^{b_i}_j\}^{r,n}_{i,j=1}.$$ Notice that $\{s^i_j\}^{r,n}_{i,j=1}$ is a set of $\C^\infty_\C(B)$--bimodule generators of the space of smooth global sections $\Gamma(E^V)$ of $E^V$.
     
   It is well--known that $\Gamma(E^V)$ is isomorphic to the space of $G$--equivariant smooth functions $$C^\infty_{\C}(P,V)^G= \{f:P\longrightarrow V \mid f \mbox{ is smooth such that } f(xA)=\alpha(A^{-1})f(x) \}$$ as $\C^\infty_\C(B)$--bimodules. Then the maps
     \begin{equation*}
     \begin{aligned}
f^i_j :P &\longrightarrow V, \\ 
x &\longrightarrow v=\sqrt{\rho_{b_i}(\pi(x))}\;\alpha(A(x)^{-1},e_j),
\end{aligned}
\end{equation*}
where $A(x)$ is the unique element of $G$ such that $\Psi_{b_i}(\pi(x),A(x))=x$, from a set of $\C^\infty_\C(B)$--bimodule generators of $C^\infty_{\C}(P,V)^G$ associated to $\{s^i_j\}^{r,n}_{i,j=1}$.  We define the smooth functions
\begin{equation*}
\begin{aligned}
f^i_{jk} :P &\longrightarrow \C \\
x &\longrightarrow \langle e_k | f^i_j(x)\rangle,
\end{aligned}
\end{equation*}
where $\langle - | -\rangle$ is the inner product that makes $\delta^V$ unitary (antilinear in the second coordinate). Now, let us consider the linear maps
$$T^i_j: V \longrightarrow C^\infty_\C(P)$$
given by $T^i_{j}(e_k)=f^i_{jk}$. A direct calculation shows that for all $x$ $\in$ $P$ and for all $A$ $\in$ $G$
$$\left(\sum^{r,n}_{l,k=1}f^{l\,\ast}_{ki}f^{l}_{kj}\right)(x)=\delta_{ij},\qquad (\Delta_P\circ T^i_j)(x,A)=((T^i_j\otimes \id_H)\circ \delta^V)(x,A).$$ Then proposition follows by taking $T^\l_1=T^1_1$, $T^\l_2=T^1_2$,..., $T^\l_n=T^1_n$, $T^\l_{n+1}=T^2_1$,..., $T^\l_{d_V}=T^r_n$ with $d_V=rn$.  
\end{proof}

It is worth mentioning that the maps $T^i_j$ agree with the {\it dualization} of $f^i_j$ via the pull--back, once the dual space $V^\#$ is identified with $V$. Moreover, to define the maps $\{T^\l_i\}^{d_V}_{i=1}$, it was necessary to have a partition of unity on the base space with smooth square roots. This is why having a $C^\ast$--algebra as the quantum base space (or one that can be 
completed to a $C^\ast$--algebra) is a sufficient condition to guarantee the existence of the maps $\{T^\l_i\}_{i=1}^{d_V}$ in the  {\it non--commutative geometrical} setting, as the reader can verify in \cite{micho3}. 

Let $\pi:P\longrightarrow B$ be a \emph{classical} principal $G$--bundle as in Proposition 4.1. Then, according to \cite{nodg}, the Gauge Principle holds: for every linear representation $\alpha^V:G\times V\longrightarrow V$, there exist a linear isomorphism $GP$ between $E^V$--valued differential forms of $B$ and basic forms of type $\alpha^V$ of $P$. This isomorphism in degree $k$
\begin{equation}
    \label{dif1}
GP^{-1}_k:\Omega^k(B)\otimes_{C^\infty_\C(B)}\Gamma(E^V)\longrightarrow \Omega^k_{\C}(P,V)^G 
\end{equation}
is given by $$ GP^{-1}_k(\mu\otimes s):=\mu\,GP^{-1}_0(s),$$  where $$(\mu\,GP^{-1}_0(s))(X_1,...,X_k):=(\pi^\#\mu)_x(X_1,...,X_k)\,GP^{-1}_0(s)(x) $$ for all $X_1,...,X_k$ $\in$ $T_xP$, with $\pi^\#\mu$ the pull--back by of $\mu$ by $\pi$, and $$GP^{-1}_0:\Gamma(E^V)\longrightarrow C^\infty_\C(P,V)^G$$ the canonical isomorphism \cite{nodg}. The existence of the map $GP$ is known as the \emph{gauge principle}. Moreover, equation (\ref{dif1}) induces the isomorphism
\begin{equation}
    \label{dif2}
GP^{'-1}:\Omega^\bullet(B)\otimes_{C^\infty_\C(B)}C^\infty_\C(P,V)^G\longrightarrow \Omega^\bullet_{\C}(P,V)^G 
\end{equation}
given by $$GP'(\mu\otimes f)=\mu\,f.$$

On the other hand, notice that the pull--back $\#$ induces an isomorphism between (for the corepresentation induced by $\alpha^V$) $$C^\infty_\C(P,V)^G\qquad \mbox{ and } \qquad \Mor(\delta^V,\Delta_P),$$ and it also induces an isomorphism between $$\Omega^\bullet_{\C}(P,V)^G\qquad\mbox{ and }\qquad \Mor(\delta^V,\Delta_\Hor),$$ where $\Delta_\Hor$ is the pull--back of the right $G$--action on the space of horizontal forms of $P$. In this way, equation (\ref{dif2}) induces the isomorphism 
\begin{equation}
    \label{dif3}
    \Upsilon^{-1}_V: \Omega^\bullet(B)\otimes_{C^\infty(B)} \Mor(\delta^V,\Delta_P)\longrightarrow \Mor(\delta^V,\Delta_\Hor)
\end{equation}
given by $$\Upsilon^{-1}_V(\mu\otimes T)=\mu\,T.$$

Let $\omega$ be a principal connection of $\pi:P\longrightarrow B$. Then, it is well--known that $\omega$ induces a linear connection on $E^V$ by means of (\cite{nodg})
\begin{equation}
    \label{dif4}
    \nabla^\omega_V:= GP\circ D^\omega \circ GP^{-1}: \Gamma(E^V)\longrightarrow \Omega^1_\C(B)\otimes_{C^\infty_\C(B)}\Gamma(E^V)
\end{equation}
and the last equation extends to the exterior derivative of $\nabla^\omega_V$:
\begin{equation}
    \label{dif5}
    d^{\nabla^\omega_V}:= GP\circ D^\omega \circ GP^{-1}: \Omega^\bullet_\C(B)\otimes_{C^\infty_\C(B)}\Gamma(E^V)\longrightarrow \Omega^\bullet_\C(B)\otimes_{C^\infty_\C(B)}\Gamma(E^V),
\end{equation}
where $D^\omega$ is the covariant derivative of $\omega$. 

Since $D^\omega$ acts on $\Omega^\bullet_{\C}(P,V)^G$ by (\cite{nodg}) 
\begin{equation}
    \label{dif6}
    D^\omega(\eta)=d\eta+\alpha^V_\ast(\omega)\wedge \eta,
\end{equation}
where $$\alpha^V_\ast(\omega)(X_x)=\left.{d\over dt}\alpha^V(\mathrm{exp}(t\omega_x(X_x))) \right|_{t=0},$$ using the pull--back $\#$ we can define an action of $D^\omega$ on $\Mor(\delta^V,\Delta_\Hor)$. Explicitly

 Therefore, we can define \emph{linear connection and its exterior derivative} on the space $$\Mor(\delta^V,\Delta_\Hor)$$ by means of
\begin{equation}
    \label{dif7}
    \nabla^\omega_V:= \Upsilon_V\circ D^\omega : \Mor(\delta^V,\Delta_P)\longrightarrow \Omega^1_\C(B)\otimes_{C^\infty_\C(B)}\Mor(\delta^V,\Delta_P)
\end{equation}
and 
\begin{equation}
    \label{dif8}
    d^{\nabla^\omega_V}:= \Upsilon_V\circ D^\omega \circ \Upsilon^{-1}_V: \Omega^\bullet_\C(B)\otimes_{C^\infty_\C(B)}\Mor(\delta^V,\Delta_P)\longrightarrow \Omega^\bullet_\C(B)\otimes_{C^\infty_\C(B)}\Mor(\delta^V,\Delta_P),
\end{equation}

\section{Quantum Principal Bundles and Dunkl Operators}

For an example of the theory presented in \cite{sald}, we will use somewhat lesser-known quantum bundles: the quantum principal $\G$--bundles developed in \cite{stheve, ds}. We shall denote them by $$\zeta=(P,B,\Delta_P).$$ These qpb's are defined by the {\it dualization} of {\it classical} principal bundles with (finite) Coxeter groups $W(H)$ as the structure group, which we shall denote them by $$\rho: P_{\mathrm{class}}\longrightarrow P_{\mathrm{class}}/W(H),$$ where $P_{\mathrm{class}}$ is the {\it classical} total space, $P_{\mathrm{class}}/W(H)$ is the {\it classical} base space, and the map $\rho:P_{\mathrm{class}}\longrightarrow P_{\mathrm{class}}/W(H)$ is the canonical projection (which is the bundle projection).

On the other hand,  horizontal forms on $\zeta$ are given by the complexification of the de--Rham graded differential $\ast$--algebra of $P_{\mathrm{class}}$, quantum differential forms of $B$ are given by the complexification of the de--Rham graded differential $\ast$--algebra of $P_{\mathrm{class}}/W(H)$, and $(\Gamma,d)$ is given by the theory of $\ast$--FODC's on finite groups \cite{ds}. It is worth mentioning that in $\zeta$ there exists a canonical qpc 
\begin{equation}
    \begin{aligned}
\omega^c:\mathfrak{qg}^\#& \longrightarrow \Omega^1(P)\\
\theta&\longmapsto \mathbbm{1}_B\otimes \theta.
\end{aligned}
\end{equation}
and its covariant derivative $D^{\omega^c}$  is exactly the de-Rham differential. Since $\omega^c$ is real and regular, $\widehat{D}^{\omega^c}=D^{\omega^c}$.

In order to give a concrete example let us focus our study on standard Dunkl connections using the development in \cite{ds}. It is worth mentioning that, since we have changed the standard definition of qpc's in order to embrace a more general theory, for us Dunkl displacements $\lambda:\mathfrak{qg}^\#\longrightarrow \Omega^1(P)$ do not need the $i=\sqrt{-1}$ factor. These qpc's are given by  $$\omega=\omega^\c+\lambda,$$ where $\lambda$ satisfies $\lambda(\pi(\phi_{\sigma_r}))=\varrho_{r}r^\#$, 
\begin{equation*}
\begin{aligned}
\varrho_{r} :P_{\mathrm{class}} &\longrightarrow \R \\
x &\longrightarrow {\kappa(r)\over \langle r|x\rangle},
\end{aligned}
\end{equation*}
$\kappa:R \longrightarrow \R$ is a multiplicative function (it is $W(H)$--invariant), $r$ $\in$ $R$ with $R$ the corresponding root system, and $r^\#$ is the element of the dual space associated with $r$. Covariant derivatives of this kind of qpc's are given by $$(D^\omega f)(x)=df(x)+\sum_{r\in R^+}\kappa(r){f(x)-f(x\sigma_r)\over\langle r|x\rangle }r$$ for every $f$ $\in$ $P$, which is a Dunkl operator in vector form \cite{stheve}. 

Let $\delta^V$ $\in$ $\T$. Then $$\nabla^\omega_V T=\sum^{n_V}_{i=1}\mu^{D^\omega T}_i\otimes_B T^V_i\,, \quad \widehat{\nabla}^\omega_V T=\sum^{n_V}_{i=1} T^V_i\otimes_B \mu^{D^\omega T}_i,$$ where $\mu^{D^\omega T}_k=\displaystyle\sum^{n_V}_{i=1}D^\omega T(e_i)f^{V\,\ast}_{ki}$ and $T(e_k)$ $\in$ $P$ \cite{sald}. Notice that Dunkl connections are multiplicative, but not regular \cite{stheve}. Furthermore, it is worth mentioning that, in differential geometry, there are no principal connections on the principal bundle $\rho: P_{\mathrm{class}}\longrightarrow P_{\mathrm{class}}/W(H)$.

On the other hand, the left and right canonical Hermitian structures are given by $$\langle T_1,T_2\rangle_\l=\sum^{n_V}_{k=1}T_1(e_k)T_2(e_k)^\ast\,,\qquad \langle T_1,T_2\rangle_\r=\sum^{n_V}_{k=1}T_1(e_k)^\ast T_2(e_k).$$

For an explicit example of Theorem $3.14$ of \cite{sald}, let us take a real Dunkl connection. These connections are characterized by $\omega=\omega^\c +\widetilde{\lambda}$, with $\widetilde{\lambda}=i\lambda$. Therefore, using equation  (\ref{generators}),
$$ \langle \nabla^\omega_{V} (T_1),T_2 \rangle_\l \,= \, \sum^{n_V}_{i=1} \mu^{D^\omega T_1}_i \langle T^V_i,T_2 \rangle_\l \,= \, \sum^{n_V}_{k=1} \left( D^\omega T_1(e_k)\right) T_2(e_k)^\ast;$$ and evaluating on any $x$ $\in$ $P_{\mathrm{class}}$ we have $$\sum^{n_V}_{k=1}  \left(dT_1(e_k)(x) T_2(e_k)^\ast(x)+i\sum_{r\in R^+}\kappa(r){T_1(e_k)(x)-T_1(e_k)(x\sigma_r)\over\langle r|x\rangle }r\, T_2(e_k)^\ast(x)\right).$$ In the same way  $$ \langle T_1, \nabla^\omega_V (T_2) \rangle_\l \,= \, \sum^{n_V}_{i=1}  \langle T_1,T_i \rangle_\l \left(\mu^{D^\omega T_2}_i\right)^\ast\,= \, \sum^{n_V}_{k=1}  T_1(e_k) \left( D^\omega T_2(e_k)\right)^\ast;$$ and evaluating on any $x$ $\in$ $P_{\mathrm{class}}$ we get $$\sum^{n_V}_{k=1}  \left(T_1(e_k)(x) dT_2(e_k)^\ast(x)-i \, T_1(e_k)(x)\sum_{r\in R^+}\kappa(r){T_2(e_k)^\ast(x)-T_2(e_k)^\ast(x\sigma_r)\over\langle r|x\rangle }r\right).$$ This implies for all $x$ $\in$ $P_{\mathrm{class}}$ that $\langle \nabla^\omega_V (T_1),T_2 \rangle_\l+ \langle T_1, \nabla^\omega_V (T_2) \rangle_\l$ is equal to
\begin{equation*}
d\langle T_1, T_2\rangle_\l(x) + i\sum_{r\in R^+} \kappa(r) {q(x)-q(x\sigma_r)\over\langle r|x\rangle }r,   
\end{equation*}
where $q=\displaystyle \sum^{n_V}_{k=1}T_1(e_k)(\sigma_r\cdot T_2(e_k))^\ast$ $\in$ $P$. However, since $T_1$, $T_2$ $\in$ $\Mor(\delta^V,\Delta_P)$, we have $\sigma_r\cdot T_i(e_j)=\displaystyle \sum^{n_V}_{k=1} \lambda^{\sigma_r}_{kj}T_i(e_k)$. Thus $q=\displaystyle \sum^{n_V}_{k=1}(\sigma_r\cdot T_1(e_k)) T_2(e_k)^\ast$ and hence $q(x)=q(x\sigma_r)$ for all $x$ $\in$ $P_{\mathrm{class}}$. In summary, we have proven explicitly that $$\langle \nabla^\omega_V (T_1),T_2 \rangle_\l+ \langle T_1, \nabla^\omega_V (T_2) \rangle_\l=d\langle T_1, T_2\rangle_\l.$$ A similar argument proves that $$\langle \widehat{\nabla}^\omega_V (T_1),T_2 \rangle_\r+ \langle T_1, \widehat{\nabla}^\omega_V (T_2) \rangle_\r=d\langle T_1, T_2\rangle_\r.$$

Defining $g^V_{ij}:W(H)\longrightarrow \C$ by $g^V_{ij}(h)=\lambda^h_{ij}$ we get

\begin{equation*}
\qtrs(g^V_{ij})=\sum^{n_V}_{k=1}  f^{V\,\ast}_{ki}\otimes_B f^V_{kj},
\end{equation*}

\noindent and by considering $\omega^c$ in equation (\ref{6.f1.6.1}), we have  for all $\theta$ $\in$ $\mathfrak{qg}^\#$  
\begin{equation}
\label{5.f2}
\qtrs(\theta)=\mathbbm{1}\otimes_{\Omega^\bullet(B)} (\mathbbm{1}\otimes \theta)-(\mathbbm{1}\otimes \theta^{(0)})\,\qtrs(\theta^{(1)}),
\end{equation}
where $\ad(\theta)=\theta^{(0)}\otimes\theta^{(1)}$.

Now let us consider a graded left $\Omega^\bullet(B)$--module isomorphism $$\F:\Omega^\bullet(B)\longrightarrow\Omega^\bullet(B)$$ which (for irreducible degree $1$ elements) is given by $$\F(\mu+x\otimes \theta)=\mu+x\lambda(\theta)+x\otimes \theta,$$ where $\lambda$ is a Dunkl displacement, $\mu$ $\in$ $\Hor^1 P$, $x$ $\in$ $P$, and $\theta$ $\in$ $\mathfrak{qg}^\#$. This map is actually a qgt $\f$ and a direct calculation shows that $\F^{\circledast}\omega^c $ is a standard Dunkl connection. By using the form of the last map $\F$, every single {\it classical} gauge transformation can be extended into a {\it quantum} one.

\section{About Adjointability on Finitely Generated Projective Modules. }

Let $(A,\cdot,\mathbbm{1},\ast)$ be a $\ast$--algebra. The map 
\begin{equation}
\label{canoleft2}
(-,-)^{d}_\l:A^{d}\times A^{d} \longrightarrow A,\qquad (\overline{x},\overline{y})\longmapsto ( \overline{x}, \overline{y})^{d}_\l:=\sum^{d}_{i=1} x_i\,y^\ast_i
\end{equation}
is the \emph{canonical} non--degenerate Hermitian structure on the free left $A$--module $$A^d=\underbrace{A\times A\times \cdots \times A}_{d-times} $$ Here, $\overline{x}=(x_1,...,x_d)$, $\overline{y}=(y_1,...,y_d)$ with $x_i$, $y_i$ $\in$ $A$.

If $M_d(A)$ denotes the space of $d\times d$ matrices with entries in $A$ and $\End(A^d)$ denotes the space of all left $A$--module endomorphisms of $A^d$, then $$\End(A^d)\cong M_d(A)$$ as rings. In fact, if $\beta:=\{\overline{e}_i \}^d_{i=1}$ is the canonical basis of $A^d$, then the previous ring isomorphism is given by $$C \longleftrightarrow C_\beta:=(C_{ij}),$$ where $C(\overline{e}_i)=\displaystyle \sum^d_{j=1} \overline{e}_j\,C_{ij}$. Furthermore, for every $\overline{x}$ $\in$ $A^d$ we have $$C(\overline{x})=\overline{x}\cdot C_\beta \qquad \mbox{ and }\qquad U(C(\overline{x}))=\overline{x}\cdot C_\beta \cdot U_\beta $$ for every $C$, $U$ $\in$ $\End(A^d)$. In addition,  for every $C_\beta$ $\in$ $M_d(A)$, a direct calculation shows that $$(\overline{x}\cdot C_\beta,\overline{y})^d_\l=(\overline{x},\overline{y}\cdot C^\dagger_\beta)^d_\l,$$ i.e., $C\cong C_\beta$ is $A$--adjointable with respect to $(-,-)^{d}_\l$. Here, $\dagger$ denotes the composition of the $\ast$ operation with the usual matrix transposition; and the $A$--adjoint operator of $C$ is the left $A$--module morphism $$C^{\star}:A^d\longrightarrow A^d$$ induced by the matrix $C^{\dagger}_\beta$.

Let  $\rho$ $\in$ $\End(A^d)$ be an idempotent element. Then $$\rho(A^d)=A^d \cdot \rho_\beta$$ is a finitely generated projective left $A$--module (\cite{lan}). Moreover, if $\rho_\beta=\rho^\dagger_\beta$, the map $(-,-)^{d}_\l$ induces a non--degenerate Hermitian structure (Section 4.3 Proposition 22 of reference \cite{lan})
$$(-,-)_\l:\rho(A^d) \times \rho(A^d) \longrightarrow A.$$ 

\begin{Proposition}
\label{a.2}
    Every left $A$--module endomorphism of $\rho(A^d)$ is $A$--adjointable  with respect to $(-,-)_\l$.
\end{Proposition}
\begin{proof}
    Let $C: \rho(A^d)\longrightarrow \rho(A^d)$ be a left $A$--module morphism and consider the map
    \begin{equation*}
        C^{\mathrm{ext}}: A^d\longrightarrow A^d,\qquad
        \overline{a} \longmapsto C(\rho(\overline{a})).
    \end{equation*}
    Since $C$ and $\rho$ are left $A$--linear maps, it follows that $C^{\mathrm{ext}}$ $\in$ $\End(A^d)$. We claim $$C=\rho\circ C^{\mathrm{ext}}\circ \rho.$$ Indeed, let $\overline{x}$ $\in$ $\rho(A^d)$. Then  $\overline{x}=\rho(\overline{a})$ for some $\overline{a}$ $\in$ $A^d$ and hence $ C(\overline{x})=C(\rho(\overline{a}))$. On the other hand,
    \begin{eqnarray*}
        (\rho\circ C^{\mathrm{ext}}\circ \rho)(\overline{x})=(\rho\circ C^{\mathrm{ext}}\circ \rho)(\rho(\overline{a}))=\rho(C^{\mathrm{ext}}(\rho(\rho(\overline{a}))))&=&\rho(C^{\mathrm{ext}}(\rho(\overline{a})))
        \\
        &=&
        \rho(C(\rho(\rho(\overline{a}))))
        \\
        &=&
        \rho(C(\rho(\overline{a}))).
    \end{eqnarray*}
   However, $\Im(C)$ $\subseteq$ $\rho(A)$, so $\rho(C(\rho(\overline{a})))=C(\rho(\overline{a}))=C(\overline{x})$ and therefore $$C(\overline{x})=(\rho\circ C^{\mathrm{ext}}\circ \rho)(\overline{x}).$$ This proves our claim. 

   In this way, by considering the ring isomorphism $\End(A^d)\cong M_d(A)$ we get $$C\cong \rho_\beta\cdot C^{\mathrm{ext}}_\beta\cdot \rho_\beta$$ and since $\rho_\beta=\rho^\dagger_\beta$, we obtain $$(\overline{a}\cdot\rho_\beta\cdot (\rho_\beta\cdot C^{\mathrm{ext}}_\beta \cdot \rho_\beta),\overline{b}\cdot\rho_\beta)_\l=(\overline{a}\cdot\rho_\beta ,\overline{b}\cdot\rho_\beta \cdot(\rho_\beta\cdot C^{\mathrm{ext}\,\dagger}_\beta \cdot \rho_\beta))_\l. $$ Hence, the induced map $$C^{\star}: \rho(A^d)\longrightarrow \rho(A^d)$$ of the matrix $\rho_\beta\cdot C^{\mathrm{ext}\,\dagger}_\beta \cdot \rho_\beta$ is the $A$--adjoint operator of $C$.
\end{proof}

Let $M$ be a finitely generated projective left $A$--module. According to Section 4.2 of reference \cite{lan}, there exists $d$ $\in$ $\N$ and an idempotent element $\rho$ $\in$ $\End(A^d)$ such that $$M\cong\rho(A^d)=A^d\cdot \rho_\beta$$ as left $A$--modules. By the last proposition, it follows that

\begin{Proposition}
    \label{a.3}
    If $\rho_\beta=\rho^\dagger_\beta$, then every left $A$--module endomorphism $C$ of $M$ is $A$--adjointable with respect to the Hermitian structure induced by $(-,-)^{d}_\l$. The $A$--adjoint operator of $C$ will be denoted by $C^{\star}$.
\end{Proposition}

Consider the free left $A$--module $A^d\otimes_A A^s$. This module is canonically isomorphic to $A^{ds}$ by means of $$\overline{\beta}_i\otimes_A \overline{\alpha}_j\;\;\longleftrightarrow\;\; \overline{e}_{(i-1)s+j},$$ where $\{\overline{e}_k\}^{ds}_{k=1}$, $\{\overline{\beta}_k\}^{d}_{k=1}$, $\{\overline{\alpha}_k\}^{s}_{k=1}$ are the canonical basis of $A^{ds}$, $A^{d}$, $A^{s}$, respectively. In this way, the canonical Hermitian structure $$(-,-)^{ds}_\l: A^{ds}\times A^{ds}\longrightarrow A^{ds}$$ defines the Hermitian structure $$(-,-)^{\otimes}_\l: (A^d\otimes_A A^s)\times (A^d\otimes_A A^s)\longrightarrow A $$ given by $$(\overline{x}\otimes_A \overline{y},\overline{u}\otimes_A \overline{v})^{\otimes}_\l:=(\overline{x}\,(\overline{y},\overline{v})^s_\l,\overline{u})^d_\l.$$

Let $M\cong \rho(A^d)=A^d\cdot \rho_\beta$, $N\cong \varsigma(A^s)=A^s\cdot \varsigma_\alpha$ be two finitely generated projective left $A$--modules with $\rho^\dagger_\beta=\rho_\beta$ and  $\varsigma^\dagger_\alpha=\varsigma$. Assume that $M\otimes_A N$ is a finitely generated projective left $A$--module too by embedding it in $A^d\otimes_A A^s\cong A^{ds}$ and its associated idempotent matrix is self--adjoint as well. Then by Proposition \ref{a.3} we get

\begin{Proposition}
    \label{a.6}
    Every left $A$--module morphism $$C: M\otimes_A N\longrightarrow M\otimes_A N$$ is $A$--adjointable with respect to the induced Hermitian structure of $(-,-)^{\otimes}_\l$ on $M\otimes_A N$.
\end{Proposition}

Of course, there are similar results for right $A$--modules and the \emph{canonical} right Hermitian structure on $A^d$ given by 
    \begin{equation}
    \label{canoright2}
     (-,-)^{d}_\r:A^{d}\times A^{d} \longrightarrow A,\qquad (\overline{x},\overline{y})\longmapsto (\overline{x},\overline{y})^{d}_\r:=\displaystyle  \sum^{d}_{i=1} x^\ast_i\,y_i.
    \end{equation}
 In particular, consider the free right $A$--module $A^d\otimes_A A^s$. This module is canonically isomorphic to $A^{ds}$ by means of $$\overline{\beta}_i\otimes_A \overline{\alpha}_j\;\;\longleftrightarrow\;\; \overline{e}_{(i-1)s+j}.$$ In this way, the canonical Hermitian structure $$(-,-)^{ds}_\r: A^{ds}\times A^{ds}\longrightarrow A^{ds}$$ defines the Hermitian structure $$(-,-)^{\otimes}_\r: (A^d\otimes_A A^s)\times (A^d\otimes_A A^s)\longrightarrow A $$ given by $$(\overline{x}\otimes_A \overline{y},\overline{u}\otimes_A \overline{v})^{\otimes}_\r:=(\overline{x},\,(\overline{y},\overline{v})^s_\r\,\overline{u})^d_\r.$$

Let $M\cong \rho(A^d)=\rho_\beta\cdot A^d $, $N\cong \varsigma(A^s)=\varsigma_\alpha \cdot A^s$ be two finitely generated projective right $A$--modules with $\rho^\dagger_\beta=\rho_\beta$ and $\varsigma^\dagger_\alpha=\varsigma$. Assume that $M\otimes_A N$ is a finitely generated projective right $A$--module too by embedding it in $A^d\otimes_A A^s\cong A^{ds}$ and its associated idempotent matrix is self--adjoint as well. Then we have

\begin{Proposition}
    \label{a.7}
    Every right $A$--module morphism $$C: M\otimes_A N\longrightarrow M\otimes_A N$$ is $A$--adjointable with respect to the induced Hermitian structure of $(-,-)^{\otimes}_\r$ on $M\otimes_A N$.
\end{Proposition}
 Of course, the proof of the last proposition is completely similar to one of Proposition \ref{a.6}.


\begin{thebibliography}{99999}
      %
  \bibitem[1]{sald}
  \textsc{Salda\~na, M, G, A.~:}\quad
  \textit{Geometry of Associated Quantum Vector Bundles and the Quantum Gauge Group,\ }
%
  \bibitem[2]{woro1}
  \textsc{Woronowicz, S, L.~:}\quad
  \textit{Compact Matrix Pseudogroups,\ }
  \textrm{Commun. Math. Phys. {\bf 111}, 613-665 (1987).}
%
 \bibitem[3]{woro2}
  \textsc{Woronowicz, S, L.~:}\quad
  \textit{Differential Calculus on Compact Matrix Pseudogroups (Quantum Groups),\ }
  \textrm{Commun. Math. Phys. {\bf 122}, 125-170 (1989).}
%
 \bibitem[4]{micho1}
  \textsc{Durdevich, M.~:}\quad
  \textit{Geometry of Quantum Principal Bundles I,\ }
  \textrm{Commun Math
Phys {\bf 175} (3), 457-521 (1996).}
%
 \bibitem[5]{micho2}
  \textsc{Durdevich, M.~:}\quad
  \textit{Geometry of Quantum Principal Bundles II,\ }
  \textrm{Rev.~Math.~Phys.~{\bf 9} (5), 531---607 (1997).}
  %
 \bibitem[6]{stheve}
  \textsc{Sontz, S, B.~:}\quad
  \textit{Principal Bundles: The Quantum Case,\ }
  \textrm{Universitext, Springer, 2015.}
%
 \bibitem[7]{war}
  \textsc{Warner, W, F.~:}\quad
  \textit{Foundations of Differentiable Manifolds and Lie Groups,\ }
  \textrm{Graduate Texts in Mathematics, Springer, 1983.}  
%
 \bibitem[8]{micho3}
  \textsc{Durdevich, M.~:}\quad
  \textit{Geometry of Quantum Principal Bundles III,\ }
  \textrm{Algebras, Groups \& Geometries.~{\bf 27}, 247--336 (2010).}
%
 \bibitem[9]{br}
  \textsc{Brzezi\'nski, T~:}\quad
  \textit{Translation Map in Quantum Principal Bundles,\ }
  \textrm{J. Geom. Phys. {\bf 20}, 349 (1996).}
 %
  \bibitem[10]{micho4}
  \textsc{Durdevich, M.~:}\quad
  \textit{Quantum Principal Bundles as Hopf-Galois Extensions,\ }
  \textrm{arXiv:q-alg/9507022v1}
%
  \bibitem[11]{micho5}
  \textsc{Durdevich, M.~:}\quad
  \textit{Quantum Gauge Transformations \& Braided Structure on Quantum Principal Bundles,\ }
  \textrm{arXiv:q-alg/9605010v1}  
 %
 \bibitem[12]{ds}
  \textsc{Durdevich, M. \&  Sontz, S, B.~:}\quad
  \textit{Dunkl Operators as Covariant Derivatives in a
Quantum Principal Bundle,\ }
  \textrm{SIGMA {\bf 9}, 040 (2013).}
%
 \bibitem[13]{dunkl}
  \textsc{Dunkl, C. F.~:}\quad
  \textit{Differential Difference Operators Associated to Reflection Groups,\ }
  \textrm{Trans. Am. Math. Soc. {\bf 3}, 457--521 (1996).}
%
%
\bibitem[15]{nodg} 
\textsc{Kol\'ar, I., Michor, P. W. \& Slovák, J.~:}
\quad
  \textit{Natural Operations in Differential Geometry,\ } internet book 

http://www.mat.univie.ac.at/~michor/kmsbookh.pdf
%
%
 \bibitem[16]{lan}
  \textsc{Landi, G.~:}\quad
  \textit{An Introduction to Noncommutative Spaces and Their Geometries,\ }
  \textrm{Springer, 1997.}  
\end{thebibliography}
\end{document}